\documentclass[11pt]{article}
\usepackage{amsmath}
\usepackage{amsthm}
\usepackage{amsfonts}
\usepackage{amssymb}
\usepackage{mathrsfs}
\usepackage{epsfig}

\parskip \smallskipamount

\newtheorem{theorem}{Theorem}[section]
\newtheorem{lemma}{Lemma}[section]

\newtheorem{proposition}{Proposition}[section]
\theoremstyle{remark}

\newtheorem{remark}{Remark}[section]
\newtheorem{example}{Example}[section]

\def\R{\mathbb{R}}
\def\Z{\mathbb{Z}}

\def\BB{\mathscr{B}}

\def\XX{\mathcal{X}}
\def\XXX{\mathscr{X}}
\def\LL{\mathcal{L}}
\def\SS{\mathcal{S}}

\renewcommand{\phi}{\varphi}
\renewcommand{\epsilon}{\varepsilon}

\newcommand{\1}{{\text{\Large $\mathfrak 1$}}}
\newcommand{\cadlag}{{c\`adl\`ag} }

\newcommand{\var}{\operatorname{var}}

\newcommand{\second}{^{\text{\tiny \rm I}}}
\newcommand{\smooth}[1]{\hat{#1}}
\newcommand{\bnd}[1]{#1^*}
\newcommand{\vr}{v}
\renewcommand{\limsup}{\varlimsup}
\renewcommand{\liminf}{\varliminf}
\newcommand{\oo}[1]{\overline{#1}}
\newcommand{\defn}{:=}

\newcommand{\keywords}[1]{ \noindent {\footnotesize
             {\small \em Keywords and phrases.} {\sc #1} } }
\newcommand{\ams}[2]{  \noindent {\footnotesize
             {\small \em AMS {\rm 2000} subject classifications.
             {\rm Primary {\sc #1}; secondary {\sc #2}} } } }

\begin{document}

\title{\bf The principle of a single big jump: discrete and continuous time
modulated random walks with heavy-tailed increments}

\author{
  \sc Serguei Foss \and \sc Takis Konstantopoulos \and \sc Stan Zachary
}
\date{
  \emph{Heriot-Watt University}\\[5mm]
  September 2005
}

\maketitle

\vspace{-5mm}

\begin{abstract}
We consider a modulated process $S$ which, 
conditional on a background process $X$,
has independent increments.
Assuming that $S$ drifts to $-\infty$ and that its increments (jumps)
are heavy-tailed (in a sense made precise in the paper), 
we exhibit natural conditions under which the asymptotics of the
tail distribution of the overall maximum of $S$
can be computed. 
We present results in discrete and in continuous time.
In particular, in the absence of modulation, the process $S$
in continuous time reduces to
a L\'evy process with heavy-tailed L\'evy measure. 
A central point of the paper is that we make full use of
the so-called ``principle of a single big jump'' in order to obtain
both upper and lower bounds. Thus, the proofs are entirely probabilistic. 
The paper is motivated by queueing and
L\'evy stochastic networks.
\\[2mm]
 \\
\ams{60G50,60G51,60K37,60F99}{60G70,60K15,90B15}
\end{abstract}

\section{Introduction}
This paper deals with the study of the asymptotic distribution of the
maximum of a random walk $S$ on the real line~$\R$, modulated by a
regenerative process, when the increments have heavy-tailed
distributions.  (By ``modulated'' we mean that, conditional on
some background process, $S$ becomes a process with independent increments.)
Our goals are 
(a) to generalise and unify existing results, 
(b) to obtain results for continuous-time modulated random walks,
and 
(c) to simplify proofs by making them entirely probabilistic, using
what we may call the {\em principle of a single big jump}, namely the
folklore fact that achieving a high value of the maximum of the random
walk is essentially due to a single very large jump. Indeed, we manage
to translate this principle into rigorous statements that make the
proofs quite transparent.

Throughout the paper, by ``tail'' we mean exclusively ``right tail'',
except where otherwise explicitly stated.  By a \emph{heavy-tailed}
distribution we mean a distribution (function) $G$ on $\R$
possessing no exponential moments: $\int_0^\infty e^{sy}G(dy)=\infty$
for all $s > 0$.  Such distributions not only abound in theory, but
are also useful in modern probabilistic modelling, in such diverse
areas as insurance risk, communication networks, and finance.

For any distribution function $G$ on $\R$, we set
$\oo{G}(y)\defn{}1-G(y)$ and denote by 
$G^{*n}$ the $n$-fold convolution of $G$ by itself. 

A distribution $G$ on $\R_+$ belongs to the class~$\SS$ of
\emph{subexponential} distributions if and only if, for all $n\ge2$,
we have $\lim_{y\to\infty}\oo{G^{*n}}(y)/\oo{G}(y)=n$.  (It is
sufficient to verify this condition in the case~$n=2$---see Chistyakov
(1964).)  This statement is easily shown to be equivalent to the
condition that, if $Y_1,\dots,Y_n$ are i.i.d.\ random variables with
common distribution $G$, then
\begin{displaymath}
  P(Y_1+\cdots+Y_n > y) \sim P(\max(Y_1,\ldots,Y_n) > y),
\end{displaymath}
a statement which already exemplifies the principle of a single big
jump. (Here, and elsewhere, for any two functions $f$, $g$ on $\R$, by
$f(y)\sim{}g(y)$ as $y\to\infty$ we mean
$\lim_{y\to\infty}f(y)/g(y)=1$; we also say that $f$ and $g$ are
\emph{tail-equivalent}.  We also write $f(y)\lesssim{}g(y)$ as
$y\to\infty$ if $\lim_{y\to\infty}f(y)/g(y)\le1$.)  The class~$\SS$
includes all the heavy-tailed distributions commonly found in
applications, in particular regularly-varying, lognormal and Weibull
distributions.

If $G_1$ and $G_2$ are distributions on $\R_+$ such that $G_1\in\SS$
and $\oo{G_2}(y)\sim{}c\oo{G_1}(y)$ as $y\to\infty$ for some
constant~$c>0$, then also $G_2\in\SS$---see Kl\"uppelberg (1988).
In particular, subexponentiality is
a tail property, a result of which we make repeated implicit use
below.  It is thus natural to extend the definition of
subexponentiality to distributions~$G$ on the entire real line~$\R$ by
defining $G\in\SS$ if and only if $G_+\in\SS$ where
$G_+(y)=G(y)\1(y\ge0)$ and $\1$ is the indicator function.  Some
further necessary results for subexponential distributions are given
in the Appendix.

We define also the class~$\LL$ of \emph{long-tailed} distributions on
$\R$ by $G\in\LL$ if and only if, for all $c$,
$\oo{G}(y+c)\sim\oo{G}(y)$ as $y\to\infty$ (it is sufficient to verify
this for any nonzero $c$).  It is known that $\SS\subset\LL$ and that
any distribution in $\LL$ is heavy-tailed---see Embrechts and Omey
(1982).  Good surveys of the basic properties of heavy-tailed
distributions, in particular long-tailed and subexponential
distributions, may be found in Embrechts \emph{et al.}\ (1997) and in
Asmussen (2000).

For any distribution~$G$ on $\R$ with finite mean, we define the
integrated (or second) tail distribution (function)~$G\second$ by
\begin{displaymath}
  \oo{G\second}(y) = 1-G\second(y)
  \defn \min\left(1,\int_y^\infty\oo{G}(z)\,dz\right).
\end{displaymath}
Note that $G\in\LL$ implies that $G\second\in\LL$, but not conversely.

Let $(\xi_n,\,n\ge1)$ be a sequence of i.i.d.\ random variables with
common distribution~$F$ on $\R$ and define the random
walk~$(S_n,\,n\ge0)$ by $S_n=\sum_{i=1}^n\xi_i$ for each $n\ge0$ (with
the convention here and elsewhere that a summation over an empty set
is zero, so that here $S_0=0$). Define $M\defn\sup_{n\ge0}S_n$.  A now
classical result (Pakes (1975), Veraverbeke (1977), Embrechts and Veraverbeke
(1982)), which we henceforth refer to as the Pakes-Veraverbeke's Theorem, states
that if $F\second\in\SS$ and if $a\defn-E\xi_1>0$ (so that in
particular $M$ is a.s.\ finite) then
\begin{equation}
  \label{eq:105}
  P(M>y) \sim \frac{1}{a}F\second(y)
  \qquad\text{as $y\to\infty$}.
\end{equation}
(Again it is the case that for most common heavy-tailed
distributions~$F$, including those examples mentioned above, we have
$F\second\in\SS$.)  The intuitive idea underlying this result is the
following: the maximum $M$ will exceed a large value $y$ if the
process follows the typical behaviour specified by the law of large
numbers, i.e.\ it's mean path, except that at some one time $n$ a jump
occurs of size greater than $y+na$; this has
probability~$\oo{F}(y+na)$; replacing the sum over all $n$ of these
probabilities by an integral yields \eqref{eq:105}; this again is the
principle of a single big jump.  See Zachary (2004) for a short proof
of \eqref{eq:105} based on this idea.

In the first part of the paper (Section~\ref{sec:discrete-time}) we
consider a sequence $(\xi_n,\,n\ge1)$ of random variables which,
conditional on another process $X=(X_n,\,n\ge1)$, are independent, and
which further are such that the conditional distribution of each
$\xi_n$ is a function of $X_n$ and otherwise independent of $n$.  We
then talk of the partial sums $S_n\defn\sum_{i=1}^n \xi_i$ as a\emph{
  modulated} random walk.  (In fact our framework includes a variety
of apparently more general processes, e.g.\ Markov additive
processes---see Remark~\ref{rem:map}.)  Our aim is to obtain the
appropriate generalisation of the result~\eqref{eq:105}.  We give
references to earlier work below.

We need to assume some asymptotic stationarity for the background
process $X$ which we take to be regenerative.  A particular case of
this is when $X$ is an ergodic Markov chain.  We also
suppose that the conditional distributions~$F_x$ given by
$F_x(y)\defn{}P(\xi_n \le y\mid{}X_n=x)$ have tails which are bounded
by that of a \emph{reference} distribution~$F$ such that
$F\second\in\SS$.  We then show (in Theorem~\ref{simple}) that, in the
case where the distributions~$F_x$ (when truncated sufficiently far
below) have means which are uniformly negative, then, under very
modest further regularity conditions, the result~\eqref{eq:105}
continues to hold with $1/a$ replaced by $C/a$.  Here $-a<0$ is now
the average (with respect to the stationary distribution of $X$) of
the above means and the constant~$C$ measures the average weight of
the tail of $F_x$ with respect to that of $F$. (The condition $-a<0$
is sufficient to ensure that $M$ is a.s.\ finite.)  In the more
general case where the distributions~$F_x$ have means of arbitrary
sign, but where $-a$ defined as above continues to be negative, we
show (in Theorem~\ref{general}) that the result~\eqref{eq:105}
continues to hold, with $1/a$ replaced by $C/a$ as above, provided
that an appropriate condition is imposed on the tail of the
distribution of the lengths of the regenerative cycles of the
process~$X$.  We give an example to show the necessity of this
condition.  Our proofs follow the probabilistic intuition of a single
big jump as defined above.  One key idea, encapsulated in a very
general result given in Section~\ref{sec:discrete-iceland} and
applicable to a wide class of processes with independent heavy-tailed
increments, is to use the result~\eqref{eq:105} of the Pakes-Veraverbeke
Theorem itself: the extremes of the increments of the general process
may be bounded by those of an unmodulated random walk, whose
increments are i.i.d.\ with negative mean; the fact that an extreme
value of the supremum of the latter process may only be obtained via a
single big jump ensures the corresponding result for the modulated
process.  Indeed we only ever use the condition~$F\second\in\SS$ in
the application of the Pakes-Veraverbeke Theorem (though we make frequent of
use the weaker condition~$F\second\in\LL$).  A preliminary version of
the discrete-time theory was given in Foss and Zachary (2002).  The
present treatment is considerably simpler and more unified, and
results are obtained under weaker conditions which are, in a sense,
demonstrated in Example~\ref{example}, optimal.

We mention several papers on the tail asymptotics of a the supremum of
a discrete-time modulated random walk with heavy-tailed increments.
Arndt (1980) considers increments with regularly varying tails
modulated by a finite-state-space Markov chain.  Alsmeyer and Sgibnev
(1999) and, independently, Jelenkovic and Lazar (1999) also consider a
finite state space Markov chain as the modulating process, and assume
that the increments of the modulated process have a subexponential
integrated tail.  Note that, for a finite Markov chain, the cycle
length distribution has an exponential tail.  Asmussen (1999)
considers a modulated random walk with an exponentially bounded
distribution of the cycle length, and assumes that both the tails and
the integrated tails of the increments of the modulated process have
subexponential distributions.  Asmussen and M{\o}ller (1999) and
Asmussen (1999) also consider a random walk with another type of
modulation.

Hansen and Jensen (2005) study the asymptotics of the maximum of a
modulated process on a finite random time horizon.  We also mention a
number of related papers on queueing systems whose dynamics may be
viewed as a kind of multi-dimensional random walk with a special type
of modulation.  Baccelli, Schlegel, and Schmidt (1999) and Huang and
Sigman (1999) consider a special type of modulation which occurs in
tandem queues and in their generalisations, and find asymptotic
results under the assumption that the tail distribution of one of the
service times strictly dominates the remainder.  A general approach to
the asymptotic study of monotone separable stochastic networks is
given by Baccelli and Foss (2004), see also Baccelli, Foss and Lelarge
(2004, 2005) for further applications.

In the second part of the paper (Section~\ref{sec:continuous-time}),
we consider the supremum of modulated continuous-time random walk,
whose jumps are similarly heavy-tailed.  The modulated continuous-time
random walk is defined as a process $(S_t,\,t\ge0)$ which, conditional
on a regenerative process $(X_t,\,t\ge0)$ has independent increments,
i.e.\ its characteristic function is given by the L\'evy-Khinchine
formula.  The parameters of the process entering the L\'evy-Khinchine
formula are therefore themselves measurable functions of the
background regenerative process.  In the absence of the background
process, $(S_t,\,t\ge0)$ becomes a L\'evy process.  Under conditions
analogous to those for the discrete-time theory, we establish similar
results for the asymptotic tail distribution of the supremum of the
process~$(S_t,\,t\ge0)$.  The continuous-time theory quite closely
parallels the discrete (and we make every attempt to keep the two
developments as similar as possible): there are, however, some
additional technicalities concerned with the ``small jumps'' and
diffusion components of the continuous-time process---these do not
contribute to the heavy-tailed distribution of the supremum; in
compensation some aspects of the theory are simpler than in discrete
time.
In particular, the proof of the lower bound in continuous-time requires 
the use of a (modulated) Poisson point process in a way that is similar
to the methods of Konstantopoulos and Richardson (2002).
  Again we require a result, given in
Section~\ref{sec:cts-iceland}, for a fairly general class of processes
with independent heavy-tailed increments.  The specialisation of this
result, under appropriate conditions, to an (unmodulated) L\'evy
process gives a simple proof of the continuous-time version of
the Pakes-Veraverbeke Theorem, different from that found in the existing
literature---see, e.g., Kl\"uppelberg, Kyprianou and Maller (2004) and
Maulik and Zwart (2005).

Some words on motivation:
heavy-tailed random variables play a significant role in the 
mathematical modeling of communication networks because the variety of
services offered by a huge system such as the Internet 
results in heterogeneous traffic. Part of the traffic concerns
small requests but other parts pose significant burden to
the system resulting in huge delays and queues.
Models of networks based on L\'evy processes--see, e.g.,
Konstantopoulos, Last and Lin (2004)
--are natural analogues
of the more-traditional Brownian networks of production and
service systems. 
To date, however, no results for the stationary distribution of
the load of stations in isolation are available. Our paper represents
a first step towards this goal. Indeed, in a L\'evy stochastic network
of feedforward type, one may see a downstream node as being
in the ``background'' of a previous node.
To apply the results of this paper to L\'evy stochastic networks is beyond
its scope and is left to a future work.

The Appendix gives some results, known and new, for the addition of
subexponential random variables, together with some other complements
to the main theory.  In particular Lemma~\ref{addtails2} extends a
well-known result for the sum of independent subexponential random
variables (Lemma~\ref{addtails}) to those which have an appropriate
conditional independence property, and is of independent interest.

\section{Modulated random walk in discrete time}
\label{sec:discrete-time}

\subsection{Introduction and main results}
\label{sec:intr-main-results}

Consider a regenerative process $X=(X_n, n\ge1)$ such that, for
each $n$, $X_n$ takes values in
some measurable space $(\XX, \XXX)$.  We say that the random walk
$(S_n,\,n\ge0)$, defined by $S_0=0$ and $S_n =\xi_1+\cdots+\xi_n$ for
$n\ge1$, is \emph{modulated} by the process $X$ if
\begin{itemize}
\item[(i)] 
conditionally on $X$, the random variables $\xi_n$, $n\ge1$, are
independent,
\item[(ii)]
for some family $(F_x,\,x\in\XX)$ of distribution functions such that,
for each $y$, $F_x(y)$ is a measurable function of $x$, we have
\begin{equation}\label{eq:95}
  P(\xi_n \le y \mid X) = P(\xi_n \le y \mid X_n) = F_{X_n}(y)
  \qquad\text{a.s.}
\end{equation}
\end{itemize}
Define
\begin{displaymath}
  M_n \defn \max(S_0, S_1, \ldots, S_n), \quad n\ge0,
  \qquad\qquad
  M \defn \sup_{n\ge 0} S_n .
\end{displaymath}
Under the conditions which we give below, $S_n \to -\infty$ a.s. as
$n\to\infty$, and the random variable~$M$ is then nondefective.  We
are interested in deriving an asymptotic expression for $P(M>y)$ as
$y\to\infty$.

\begin{remark}\label{rem:map}
  In fact nothing below changes if, in \eqref{eq:95}, we allow the
  distribution of $\xi_n$ to depend on the history of the modulating
  process~$X$ between the last regeneration instant prior to time~$n$
  and the time~$n$ itself.  This possible relaxation can either be
  checked directly, or brought within the current structure by
  suitably redefining the process~$X$.  Thus in particular our
  framework includes Markov additive processes.
\end{remark}

The regeneration epochs of the modulating process~$X$ are denoted by
$0\le{}T_0<T_1<\dots$. By definition, the cycles
$((X_n, T_{k-1} < n \le T_{k}),\, k\ge1)$ are i.i.d.\ and independent
of the initial cycle $(X_n,\, 0 < n \le T_0)$.  Define also
$\tau_0\defn{}T_0$ and $\tau_k\defn{}T_k-T_{k-1}$ for $k\ge1$, so that
$(\tau_k,\,k\ge0)$ are independent and $(\tau_k,\,k\ge1)$ are
identically distributed.  Assume that $E\tau_1<\infty$.
For each $n\ge0$, let $\pi_n$ be the distribution of
$X_n$, and define, as usual, the stationary probability measure
\begin{displaymath}
  \pi(B) \defn \frac{E\sum_{n=T_0+1}^{T_1} \1(X_n \in B)}{E\tau_1},\qquad
  B \in \BB({\mathcal X}).
\end{displaymath}


Each distribution $F_x$, $x\in\XX$, will be assumed to have a
finite mean
\begin{equation}
  \label{eq:1}
  a_x \defn \int_\R y F_x(dy).
\end{equation}
The family of such distributions will be assumed to satisfy the
following additional conditions with respect to some \emph{reference}
distribution~$F$ with finite mean and some measurable function
$c:\XX\to[0,1]$:
\begin{flalign*}
  \text{\bf (D1)} & \quad
  \oo{F_x}(y)  \le \oo{F}(y),
  \quad \text{ for all } y \in \R, \qquad x\in\XX, &\\
  \text{\bf (D2)} & \quad
  \oo{F\second_x}(y)  \sim c(x) \oo{F\second}(y)
  \quad \text{ as } y \to \infty, \qquad x\in\XX, &\\
  \text{\bf (D3)} & \quad
  \text{$a: = -\int_\XX a_x \pi(dx)$ is finite and strictly positive.}
\end{flalign*}

\begin{remark}\label{rem:gen}
  The condition (D1) is no less restrictive than the condition
  \[
  \limsup_{y\to\infty} \sup_{x\in\XX} \frac{\oo{F_x}(y)}{\oo{F}(y)} < \infty,
  \]
  in which case it is straightforward to redefine $F$, and then $c$,
  so that (D1) and (D2) hold as above.
\end{remark}
\begin{remark}\label{rem:simple}
  A sufficient condition for (D2) to hold is that, for all $x\in\XX$,
  we have $\oo{F_x}(y)\sim c(x)\oo{F}(y)$ as $y\to\infty$.  However,
  in order to obtain our main results we shall require
  Lemma~\ref{upper} below to be established under the weaker
  condition~(D2) as stated above.  (The proof of Theorem~\ref{general}
  utilises the fact that, when $F\second\in\LL$ as required there, the
  condition~(D2) is preserved when any of the distributions~$F_x$ is
  shifted by a constant. This is not true if (D2) is replaced by the
  strengthened version above, unless we further assume $F\in\LL$---an
  assumption which we do not wish to make!)
\end{remark}
\begin{remark}\label{rem:sign}
  We impose no \emph{a priori} restrictions on the signs of the $a_x$,
  other than that given by the condition~(D3).  The latter condition
  is trivially satisfied in the case where all the $a_x$ are strictly
  negative.  (The introduction of the minus sign in the definition of
  $a$ is for convenience in the statement of our results.)
\end{remark}
It follows from the regenerative structure of $X$ and from (D1)
and (D3), that
\begin{equation}\label{eq:97}
  S_n/n \to -a \quad \text{ as } n \to \infty, \quad \text{ a.s.}
\end{equation}
(See the Appendix for a proof of this result.)

Thus, in particular, $S_n\to-\infty$ as $n\to\infty$ and $M$ is
nondefective as required.

For each $x\in\XX$ and $\beta>0$, define
\begin{equation}
  \label{eq:90}
  a_x^\beta \defn a_x - \int_{-\infty}^{-\beta} (y+\beta)\,F_x(dy)
  = \int_\R (y \vee -\beta)\,F_x(dy);
\end{equation}
note that $a_x^\beta\ge{}a_x$.  Define also
\begin{equation}
  \label{eq:2}
  \kappa\defn \lim_{\beta \to \infty} \sup_{x\in\XX}a_x^\beta.
\end{equation}
Note that, from \eqref{eq:1} and the condition~(D1), $\kappa$ is a
real number between $-a$ and $\mu$, where $\mu$ is the mean of the
reference distribution~$F$.  In the case where the distributions
$F_x$, $x\in\XX$, satisfy the uniform integrability condition
\begin{displaymath}
  \lim_{\beta \to \infty} \sup_{x\in\XX}
  \int_{-\infty}^{-\beta} |y|~F_x(dy) = 0,
\end{displaymath}
we have $\kappa=\sup_{x\in\XX}a_x$.

Define also $C\in[0,1]$ by
\begin{equation}
  \label{eq:3}
  C \defn \int_\XX c(x) \pi(dx).
\end{equation}
Theorem~\ref{simple} below gives our main result in the
case~$\kappa<0$.
\begin{theorem}
  \label{simple}
  Suppose that (D1)--(D3) hold, that $F\second\in\SS$, and
  that $\kappa<0$.  Then
\begin{displaymath}
  \lim_{y\to\infty}\frac{P(M>y)}{\oo{F\second}(y)}
  = \frac{C}{a}.
\end{displaymath}
\end{theorem}
In order to extend Theorem~\ref{simple} to the case where the sign
of $\kappa$ may be arbitrary, we require an additional condition
regarding the (tail) distributions of the lengths of the regenerative
cycles.  The condition we need is:

\text{\bf (D4)} \quad For some nonnegative $b>\kappa$,
\begin{equation}
  \label{eq:4}
  P(b\tau_0 > n)  = o(\oo{F\second}(n)), \quad P(b\tau_1>n) = o(\oo{F}(n)),
\qquad \text{as $n\to\infty$.}
\end{equation}
Note that if \eqref{eq:4} is satisfied for some nonnegative $b$, then
it is also satisfied for any smaller value of $b$.  In the case
$\kappa<0$ the condition~(D4) is always trivially satisfied by taking
$b=0$.  Hence Theorem~\ref{simple} is actually a special case of the
general result given by Theorem~\ref{general} below.

\begin{theorem}
\label{general}
Suppose that (D1)--(D4) hold and that $F\second\in\SS$.  Then
\[
\lim_{y\to\infty} \frac{P(M>y)}{\oo{F\second}(y)} = \frac{C}{a}.
\]
\end{theorem}

In Section~\ref{sec:example} we give an example to show the necessity
of the assumption (D4).

\subsection{A uniform upper bound for discrete-time processes with
  independent increments}
\label{sec:discrete-iceland}

Our proofs require several uses of the following proposition, which is
new and may be of independent interest.  This,
under appropriate conditions, provides an upper bound for the
distribution of the supremum of a random walk with independent
increments. This bound is not simply asymptotic and further has an
important uniformity property.  No regenerative structure is assumed,
and the result is therefore of independent interest.

\begin{proposition}
  \label{iceland}
  Let $F$ be a distribution function on $\R$ such that
  $\int_0^\infty\oo{F}(y)\,dy<\infty$ and whose integrated tail
  $F\second\in\SS$.  Let $\alpha$, $\beta$ be given positive real
  numbers.  Consider any sequence~$(\xi_n,\,n\ge1)$ of independent
  random variables such that, for each $n$, the distribution~$F_n$ of
  $\xi_n$ satisfies the conditions
  \begin{align}
    \oo{F_n}(y) & \le \oo{F}(y)
    \qquad\text{for all $y\in\R$,} \label{eq:5}\\
    \int_\R (z \vee -\beta)~F_n(dz) & \le - \alpha. \label{eq:6}
  \end{align}
  Let $M \defn \sup_{n\ge0}\sum_{i=1}^n\xi_i$.  Then there exists a
  constant $r$ depending on $F$, $\alpha$ and $\beta$ only, such that,
  for all sequences~$(\xi_n,\,n\ge1)$ as above,
  \begin{equation}
    \label{unif}
    P(M > y) \le r \oo{F\second}(y)
    \qquad\text{ for all $y$}.
  \end{equation}
\end{proposition}

\begin{proof}
  Consider any sequence $(\xi_n,\,n\ge1)$ as above.  We assume,
  without loss of generality, that $\xi_n \ge -\beta$, a.s.\ for all
  $n$ (for, otherwise, we can replace each $\xi_n$ by
  $\max(\xi_n,-\beta)$).  We now use a coupling construction.  Let
  $(U_n,\,n\ge1)$ be a sequence of i.i.d.\ random variables with
  uniform distribution on the unit interval $(0,1)$.  For each $n$,
  let $F_n^{-1}(y)=\sup\{z:F_n(z)\le y\}$ be the generalised inverse
  of $F_n$, and define similarly $F^{-1}$.  Let
  \[
  \xi_n \defn F_n^{-1}(U_n), \quad \eta_n \defn F^{-1}(U_n).
  \]
  Then $\xi_n$ has distribution $F_n$, $\eta_n$ has distribution $F$
  and $\xi_n\le\eta_n$ a.s.  Choose a constant $y^*$ sufficiently
  large, such that
  \begin{equation}\label{eq:7}
    m \defn E[\eta_1\1(\eta_1>y^*)] \le \alpha/4
    \qquad
    \text{and}
    \qquad
    \max(1,\beta) P(\eta_1 > y^*) \le \alpha/4.
  \end{equation}
  Let $\epsilon = P(\eta_1>y^*)$ and let $K_0 = m/\epsilon + 1$.  For
  each $n$, define the random variables
  \begin{align}
    \delta_n & \defn \1(\eta_n>y^*) \label{eq:8}\\
    \phi_n & \defn \xi_n(1 - \delta_n) + K_0\delta_n \label{eq:9}\\
    \psi_n & \defn (\eta_n-K_0)\delta_n. \label{eq:10}
  \end{align}
  Note that, from \eqref{eq:6}, \eqref{eq:7}--\eqref{eq:10}, and our
  assumption that
  $\xi_n \ge -\beta$, a.s.,
  \begin{equation}
    \label{eq:11}
    E\phi_n \le E\xi_n + (\beta+K_0)E\delta_n
    \le -\alpha + (\beta+1)\epsilon + m
    \le -\alpha/4
  \end{equation}
  and
  \begin{equation}
    \label{eq:12}
    E\psi_n = m - K_0\epsilon = -\epsilon < 0.
  \end{equation}
  Note also that $(\delta_n,\,n\ge1)$ and $(\psi_n,\,n\ge1)$ are
  both sequences of i.i.d.\ random variables.  For each $n\ge0$, define
  $S^\phi_n\defn\sum_{i=1}^n\phi_i$ and
  similarly $S^\psi_n\defn\sum_{i=1}^n\psi_i$.
  Define also $M^\phi\defn\sup_{n\ge0}S^\phi_n$ and
  $M^\psi\defn\sup_{n\ge0}S^\psi_n$.  (It will follow below that
  $M^\phi$ and $M^\psi$ are almost surely finite.)  From \eqref{eq:9},
  \eqref{eq:10}, and since $\xi_n\le\eta_n$ a.s., it follows that, for
  each $n$, $\xi_n\le\phi_n+\psi_n$, and so
  \begin{equation} \label{eq:13} M \le \sup_{n\ge0}(S^\phi_n +
    S^\psi_n) \le M^\phi + M^\psi.
  \end{equation}
  Given any realisation of the two sequences $(\delta_n,\,n\ge1)$ and
  $(\phi_n,\,n\ge1)$ such that $\sum_n\delta_n = \infty$, the
  conditional distribution of $M^\psi$ coincides with that of the
  supremum of the partial sums of an i.i.d.\ sequence
  $(\psi'_n,\,n\ge1)$ where
  \begin{equation}\label{eq:98}
    P(\psi'_n\in\cdot) \defn P(\eta_1-K_0\in\cdot\,\mid\,\delta_1=1).
  \end{equation}
  It follows from (\ref{eq:12}) that $E\psi'_1=-1$.  Since
  $\sum_n\delta_n=\infty$ a.s., it follows that the random variable
  $M^\psi$ is finite a.s.\ and does not depend on the joint
  distribution of the random variables $(\delta_n,\varphi_n,\ n\ge1)$.
  In particular, $M^{\psi}$ and $M^\phi$ are independent random
  variables.  Further, since $F\second\in\SS\subset\LL$,
  it follows from \eqref{eq:98} that the common
  distribution~$F_\psi$ of the random variables~$\psi'_n$ satisfies
  $\oo{F_\psi\second}(y)\sim\oo{F\second}(y)/\epsilon$ as
  $y\to\infty$.  Hence, by the Pakes-Veraverbeke Theorem,
  \begin{equation}
    \label{eq:14}
    P(M^\psi > y) \sim \frac{1}{\epsilon} \oo{F\second}(y)
    \qquad
    \text{as $y\to\infty$.}
  \end{equation}
  We now consider the tail distribution of $M^\phi$ and show that this
  is exponentially bounded.  For each $n$, let $F^\phi_n$ be the
  distribution of $\phi_n$.  We show first how to choose a
  constant~$s$, depending on $F, \alpha, \beta$ only, such that the
  process $\exp sS^\phi_n$ is a supermartingale.  For this we require
  that, for all $n$,
  \begin{equation}
    \label{eq:15}
    \frac{1}{s}\int_{-\infty}^\infty\left(e^{sz}-1\right)F^\phi_n(dz) \le 0.
  \end{equation} 
  From \eqref{eq:8}, \eqref{eq:9}, and our assumption that
  $\phi_n\ge-\beta$ a.s., it follows that, for all $n$,
  \begin{displaymath}
    |\phi_n| \le K\defn \max(\beta, y^*, K_0).
  \end{displaymath}
  From this, and the inequality $e^{sz}\le1+sz+s^2K^2e^{sK}$, valid
  for any $s\ge0$ and for any $z$ such that $|z|\le{}K$, it follows
  that the left side of \eqref{eq:15} is bounded above by
  $E\phi_n+sK^2e^{sK}$, which, by \eqref{eq:11}, is less than or equal
  to zero for any $s>0$ such that $sK^2e^{sK}\le-\alpha/4$.

  Thus we fix such an $s$, depending only on $F, \alpha, \beta$ as
  required.  It now follows by the usual argument involving the
  martingale maximal inequality that, for $y\ge0$,
  \begin{equation}
    \label{eq:16}
    P(M^{\phi}>y) \le e^{-sy}.
  \end{equation}
  Let $\zeta$ be a random variable which has tail distribution
  $e^{-sy}$ and which is independent of everything else.  Since
  $M^\phi$ and $M^\psi$ are independent, it follows from \eqref{eq:13}
  that
  \begin{equation}\label{eq:17}
    P(M > y) \le P(M^{\psi} + \zeta > y).
  \end{equation}
  Further, from Lemma~\ref{addtails},
  \begin{displaymath}
    \lim_{y\to\infty}
    \frac{P(M^\psi+\zeta > y)}{\oo{F\second}(y)}
    =
    \lim_{y\to\infty}
    \frac{P(M^{\psi}>y)}{\oo{F\second}(y)} = \frac{1}{\epsilon}.
  \end{displaymath}
  and so there exists $r$ such that, for all $y\ge0$,
  \begin{equation}\label{eq:18}
    P(M^\psi+\zeta > y) \le r\oo{F\second}(y).
  \end{equation}
  Finally, note that the distributions of $M^\psi$ and $\zeta$, and so
  also the constant~$r$, depend on $F$, $\alpha$ and $\beta$ only, so
  that the required result now follows from \eqref{eq:17} and
  \eqref{eq:18}.
\end{proof}

\subsection{Proofs of Theorems~\ref{simple} and \ref{general}}
\label{sec:proofs-main-results}

We now return to the model and assumptions of
Section~\ref{sec:intr-main-results} and prove our main results.

We give first the following technical lemma, which will be required
subsequently.
\begin{lemma}
  \label{auxlem}
  Suppose that $F\second\in\LL$ and that $d_1$, $d_2$ are constants
  such that $d_2>0$.  Then
  \begin{equation} \label{eq:19}
    \sum_{n\ge1} \oo{F}(y+d_1+d_2n)
    \sim \frac{1}{d_2}  \oo{F\second}(y)
    \qquad\text{as $y\to\infty$,}
  \end{equation}
  and, for any real sequence $\delta_n$ such that
  $\delta_n\downarrow0$,
  \begin{equation} \label{eq:20}
    \sum_{n\ge1} \delta_n \oo{F}(y+d_1+d_2n)
    = o\left(\oo{F\second}(y)\right)
    \qquad\text{as $y\to\infty$.}
  \end{equation}
  The conditions~(D1) and (D2) further imply that
  \begin{equation} \label{eq:21}
    \sum_{n\ge1} \int_\XX\pi_n(dx) \oo{F_x}(y+d_1+d_2n)
    \sim \frac{C}{d_2} \oo{F\second}(y)
    \qquad\text{as $y\to\infty$.}
  \end{equation}
\end{lemma}

\begin{proof}
  The results \eqref{eq:19} and \eqref{eq:20} are elementary
  consequences of the condition $F\second\in\LL$, and, in each case,
  the approximation of a sum by an integral.  Detailed proofs may be
  found in Foss and Zachary (2002).  We prove \eqref{eq:21} under the
  assumption that the regenerative process~$X$ is aperiodic, so that
  the distance $||\pi_n-\pi||$ between $\pi_n$ and $\pi$ in the total
  variation norm tends to zero---the modifications required to deal
  with the periodic case are routine.  Then
  \begin{equation}
    \label{eq:22}
    \sum_{n\ge1} \int_\XX \pi_n(dx) \oo{F_x}(y+d_1+d_2n)
    \sim
    \sum_{n\ge1} \int_\XX \pi(dx) \oo{F_x}(y+d_1+d_2n) 
    \qquad\text{as $y\to\infty$,}
  \end{equation}
  since the absolute value of the difference between the left and
  right sides of \eqref{eq:22} is bounded by
  $\sum_n||\pi_n-\pi||\oo{F}(y+d_1+d_2n)$, which, by \eqref{eq:20}, is
  $o\left(\oo{F\second}(y)\right)$ as $y\to\infty$.  Further, using
  the condition~(D1), for $y$ sufficiently large that
  $\oo{F\second}(y+d_1)<1$, the right side of \eqref{eq:22} is bounded
  above and below by
  \begin{displaymath}
    \frac{1}{d_2} \int_\XX \pi(dx) \oo{F\second_x}(y+d_1)
    \quad \text{and} \quad
    \frac{1}{d_2} \int_\XX \pi(dx) \oo{F\second_x}(y+d_1+d_2)
  \end{displaymath}
  respectively.  From the conditions~(D1), (D2) and the dominated
  convergence theorem, for any constant~$d$,
  \begin{displaymath}
    \int_\XX \pi(dx) \oo{F\second_x}(y+d)
    \sim
    \oo{F\second}(y+d) \int_\XX \pi(dx) c(x)
    \qquad\text{as $y\to\infty$.}
  \end{displaymath}
  The result \eqref{eq:21} now follows from the condition
  $F\second\in\LL$
\end{proof}


The following lemma gives an asymptotic lower bound for $P(M>y)$.
This result is also proved in Foss and Zachary (2002), but we give
here for completeness a short, simplified proof---see also Zachary
(2004).

\begin{lemma}
  \label{lower}
  Suppose that (D1)--(D3) hold and that $F\second\in\LL$.  Then
  \[
  \liminf_{y\to\infty} \frac{P(M>y)}{\oo{F\second}(y)} \ge 
  \frac{C}{a}.
  \]
\end{lemma}

\begin{proof}
  Given $\epsilon>0$, by the weak law of large numbers we may choose a
  constant~$l_0$ sufficiently large that if, for each~$n$, we define
  $l_n=l_0+(a+\epsilon)n$, then
  \begin{equation}
    \label{eq:99}
    P(S_n>-l_n) > 1-\epsilon.
  \end{equation}
  For any fixed $y\ge0$ and each $n\ge1$, define
  \begin{math}
    A_n\defn{}\{M_{n-1}\le{}y,\,S_{n-1}>-l_{n-1},\,\xi_n>y+l_{n-1}\}.
  \end{math}
  Since, conditional on the background process~$X$, the random
  variables~$\xi_n$ are independent, it follows that
  \begin{align}
    P(A_n)
    & = E[\1(\{M_{n-1}\le y,\,S_{n-1}>-l_{n-1}\})\oo{F_{X_n}}(y+l_{n-1})]
    \nonumber\\
    & \ge E[\oo{F_{X_n}}(y+l_{n-1})]
    - P(\{M_{n-1} > y\}\cup\{S_{n-1}\le-l_{n-1}\})\oo{F}(y+l_{n-1})
    \label{eq:101}\\
    & \ge E[\oo{F_{X_n}}(y+l_{n-1})]
    - [P(M > y) + \epsilon]\oo{F}(y+l_{n-1}), \label{eq:102}
  \end{align}
  where \eqref{eq:101} follows from the condition~(D1) and
  \eqref{eq:102} follows from \eqref{eq:99}.
  Since also the events $A_n$, $n\ge1$, are disjoint and each is
  contained in the event~$\{M>y\}$, it follows that
  \begin{align}
    P(M>y)
    & \ge \sum_{n\ge1} E[\oo{F_{X_n}}(y+l_{n-1})]
    - [P(M > y) + \epsilon]\sum_{n\ge1}\oo{F}(y+l_{n-1})
    \nonumber\\
    & = (1+o(1))\frac{C}{a+\epsilon}\oo{F\second}(y)
    - (1+o(1))[P(M > y) + \epsilon]%
    \frac{\oo{F\second}(y)}{a+\epsilon} \label{eq:103}\\
    & = (1+o(1))\frac{C-\epsilon}{a+\epsilon}\oo{F\second}(y),
    \label{eq:104}
  \end{align}
  as $y\to\infty$, where \eqref{eq:103} follows from Lemma
  \ref{auxlem}, and \eqref{eq:104} follows since $P(M > y)\to0$ as
  $y\to\infty$.  The required result now follows by letting $\epsilon$
  tend to zero.
\end{proof}

\begin{remark}
  As in the Pakes-Veraverbeke Theorem for unmodulated random walks, the
  intuitive idea underlying the above result is the following: the
  maximum $M$ will exceed a large value $y$ if the process follows the
  typical behaviour specified by the law of large numbers, i.e.\ it's
  mean path, except that at any time $n$ a jump occurs of size greater
  than $y+na$; this has probability~$E[\oo{F_{X_n}}(y+na)]$, and so
  the bound is now given by the use of \eqref{eq:21}.  

  We shall argue similarly for the upper bound: if $M$ exceeds a large
  value $y$ then it \emph{must} be the case that a single jump exceeds
  $y$ plus the typical behaviour of the process.  We now proceed to
  making this heuristic more precise.
\end{remark}

We consider first, in Lemma~\ref{upper} below, the upper bound for the
relatively simple case $\kappa<0$.  This result may be combined with
the lower bound of Lemma~\ref{lower} to give the exact asymptotics in
this case (Theorem~\ref{simple}).  We then use the \emph{result}
of Lemma~\ref{upper} to extend the upper bound, in the proof of
Theorem~\ref{general}, to general $\kappa$, thereby obtaining the
exact asymptotics in this case also.  

\begin{lemma}
  \label{upper}
  Suppose that (D1)--(D3) hold, that $F\second\in\SS$, and that
  $\kappa<0$.  Then
  \begin{displaymath}
    \limsup_{y\to\infty}\frac{P(M>y)}{\oo{F\second}(y)}
    \le \frac{C}{a}.
  \end{displaymath}
\end{lemma}

\begin{proof}
  For given (small) $\epsilon>0$, and (large) $u_0>0$, for
  each~$n\ge0$ define $u_n = u_0-(a-\epsilon)n$.  Define the stopping
  time
  \begin{displaymath}
    \sigma = \inf\{n\ge 0:~ S_n > u_n\}.
  \end{displaymath}
  Since $S_n/n\to-a$ a.s., it follows that (for fixed $\epsilon$)
  \begin{equation}
    \label{eq:24}
    P(\sigma<\infty) \to 0 \qquad\text{as $u_0\to\infty$}.
  \end{equation}
  Note that $S_\sigma$ and $M_\sigma=\max_{0\le{}n\le{}\sigma}S_n$ are
  only defined on $\{\sigma<\infty\}$.  Here, and elsewhere, we use
  the convention that any probability of an event involving random
  variables such as $S_\sigma$ or $M_\sigma$ is actually the
  probability of the same event intersected by $\{\sigma<\infty\}$,
  e.g.\ $P(M_\sigma>y)\defn{}P(M_\sigma>y,\,\sigma<\infty)$.

  Since $S_n \le u_n$ for all $n < \sigma$, we have, for $y>u_0$,
  \[
  \{M_\sigma>y\} = \{S_\sigma > y\} = \{S_{\sigma-1} + \xi_\sigma >
  y\} \subseteq \{\xi_\sigma >y-u_{\sigma-1)}\},
  \]
  and hence
  \begin{align}
    P(M_\sigma>y) = P(S_\sigma > y) & \le \sum_{n=1}^\infty
    P(\xi_n >y-u_{n-1}) \nonumber\\
    &= \sum_{n=1}^\infty \int_{\XX} \pi_n(dx) \oo{F_x}(y-u_{n-1} \nonumber)\\
    & \sim \frac{C}{a-\epsilon} \oo{F\second}(y), \label{eq:25}
  \end{align}
  as $y\to\infty$, where the last equivalence follows from
  Lemma~\ref{auxlem}.

  Since $\kappa<0$ it follows from \eqref{eq:2} that we can choose any
  $\alpha\in(0,-\kappa)$ and then $\beta>0$ sufficiently large that
  \begin{equation}\label{eq:26}
    \sup_{x\in\XX}\int_\R (z \vee -\beta)~F_x(dz) \le - \alpha.
  \end{equation}
  On the set $\{\sigma<\infty\}$ define the sequence of random
  variables $(\xi^\sigma_n\,\,n\ge1)$ by
  $\xi^\sigma_n=\xi_{\sigma+n}$; let
  $M^\sigma=\sup_{n\ge0}\sum_{i=1}^n\xi^\sigma_i$.  Conditional on the
  background process~$X$ and any finite value of $\sigma$, the
  sequence $(\xi^\sigma_n,\,n\ge1)$ consists of independent random
  variables which, from~\eqref{eq:26}, satisfy the conditions
  \eqref{eq:5} and \eqref{eq:6} of Proposition~~\ref{iceland} (with
  $F$, $\alpha$ and $\beta$ as defined here).  It therefore follows
  from that proposition that there exists $r$, depending on $F$,
  $\alpha$ and $\beta$ only, such that, for all $x$, all finite $n$
  and all $y\ge0$,
  \begin{equation}
    \label{eq:27}
    P(M^\sigma>y \mid X=x,\,\sigma=n) \le r \oo{F\second}(y);
  \end{equation}
  further, conditional on $X=x$ and $\sigma=n$, the random variables
  $M_\sigma$ and $M^\sigma$ are independent.

  Let $\widetilde{M}$ be a random variable, \emph{independent of all
    else}, with tail distribution
  \begin{equation}
    \label{eq:28}
    P(\widetilde{M}>y) = 1 \wedge r \oo{F\second}(y),
  \end{equation}
  Observe that, for $y>u_0$, we have
  $M=S_\sigma+M^\sigma=M_\sigma+M^\sigma$.  By conditioning on $X$ and
  each finite value of $\sigma$, it follows from (\ref{eq:27}),
  (\ref{eq:28}) and the above conditional independence that
  \begin{align}
    P(M > y)
    & = P(M_\sigma+M^\sigma > y,\, \sigma<\infty)\nonumber\\
    & \le P(M_\sigma+\widetilde{M} > y,\, \sigma<\infty)\nonumber\\
    & = P(\sigma < \infty)P(M_\sigma+\widetilde{M} > y \,\mid\,
    \sigma<\infty)
    \label{eq:29}
  \end{align}
  Also, from \eqref{eq:25},
  \begin{equation}\label{eq:30}
    \limsup_{y\to\infty}%
    \frac{P(M_{\sigma}>y \,|\,\sigma<\infty)}{\oo{F\second}(y)}
    \le \frac{C}{(a-\epsilon) P(\sigma < \infty)}.
  \end{equation}
  From \eqref{eq:28}, \eqref{eq:29}, \eqref{eq:30}, the independence
  of $\widetilde{M}$ from all else, and Lemma~\ref{addtails},
  \begin{align}
    \limsup_{y\to\infty}\frac{P(M>y)}{\oo{F\second}(y)} & \le
    P(\sigma<\infty) \left( \frac{C}{(a-\epsilon) P(\sigma < \infty)}
      + r
    \right) \nonumber\\
    & = \frac{C}{a-\epsilon} + rP(\sigma<\infty). \label{eq:31}
  \end{align}
  It follows from \eqref{eq:24} that as $u_0\to\infty$ the second term
  on the right side of \eqref{eq:31} tends to $0$.  The required
  result now follows on letting also $\epsilon\to0$.
\end{proof}
 
\begin{proof}[Proof of Theorem~\ref{simple}]
  This is now immediate from Lemmas~\ref{lower} and \ref{upper}.
\end{proof}

\begin{proof}[Proof of Theorem~\ref{general}]
  Let nonnegative $b>\kappa$ be such that the condition~(D4) holds.
  Choose
  \begin{equation}
    \label{eq:32}
    \delta\in(0,\min(a,b-\kappa))
  \end{equation}
  and choose $\epsilon\in(0,a-\delta)$.  Note that, from the
  condition~(D3),
  \begin{displaymath}
    \int_\XX(a_x + \delta)\pi(dx)= -a + \delta < -\epsilon.
  \end{displaymath}
  It now follows from the from the definition~\eqref{eq:2} of
  $\kappa$, and since also $b>\kappa$, that we may choose $\beta>0$
  sufficiently large that
  \begin{gather}
    \int_\XX(a_x^\beta + \delta)\pi(dx) < -\epsilon, \label{eq:92}\\
    a_x^\beta + \delta \le b \qquad\text{for all $x\in\XX$},
    \label{eq:93}
  \end{gather}
  (where $a_x^\beta$ is as defined by \eqref{eq:90}).  Let $b_x$ be a
  measurable function on $\XX$ such that, for some sufficiently large
  $d>0$,
  \begin{gather}
    \max(-d,\, a_x^\beta+\delta) \le b_x \le b,
    \qquad x\in\XX, \label{eq:33}\\
    \int_\XX b_x \pi(dx) = -\epsilon. \label{eq:34}
  \end{gather}
  (To see that such a function $b_x$ exists, note that,
  from~\eqref{eq:92}, we may choose $d$ sufficiently large that
  \begin{displaymath}
    \int_\XX\max(-d,\, a_x^\beta+\delta) \pi(dx) < -\epsilon;
  \end{displaymath}
  since also $b>-\epsilon$, we may, for example, satisfy \eqref{eq:33}
  and \eqref{eq:34} by choosing $b_x=\max(s,-d,\,a_x^\beta+\delta)$
  for the appropriate constant $s\in(-d,b)$.)

  Define, for each $n\ge1$,
  \begin{equation}
    \label{eq:35}
    \smooth{\xi}_n = \xi_n - b_{X_n}.
  \end{equation}
  Note that, conditional on the modulating process~$X$, the random
  variable~$\smooth{\xi}_n$ has distribution function
  $\smooth{F}_{X_n}$ where, for each $x\in\XX$ and each $y\in\R$,
  $\smooth{F}_x(y)=F_x(y+b_x)$.  Since also
  $F\second\in\SS\subset\LL$,
  the family of distributions
  $(\smooth{F}_x,\,x\in\XX)$ satisfies the conditions~(D1) and (D2)
  with $F$ replaced by $\smooth{F}$ where $\smooth{F}(y)=F(y-d)$.

  Since also, for each $x\in\XX$, the distribution $\smooth{F}_x$ has
  mean $\smooth{a}_x\defn{}a_x-b_x$ and $b_x\le{}b$, it follows
  that, for $\beta'=b+\beta$,
  \begin{equation}
    \label{eq:36}
    \int_\R (z \vee -\beta')\smooth{F}_x(dz) \le a_x^\beta - b_x \le -\delta,
  \end{equation}
  where the second inequality above also follows from~\eqref{eq:33}.
  Lastly, it follows from the condition~(D3) and \eqref{eq:34} that,
  for each $x\in\XX$,
  \begin{displaymath}
    \int_\XX \smooth{a}_x \pi(dx) = -a + \epsilon < 0.
  \end{displaymath}
  The process~$(\smooth{S}_n,\,n\ge0)$ given by 
  $\smooth{S}_n=\sum_{i=1}^n\smooth{\xi}_i$, for each $n\ge0$, thus
  satisfies all the conditions associated with Lemma~\ref{upper},
  where $F$ is replaced by $\smooth{F}$, $\kappa$ is replaced by the
  appropriate $\smooth{\kappa}$ with, from~\eqref{eq:36},
  $\smooth{\kappa}\le-\delta$, and $a$ is replaced by $a-\epsilon$.
  Since also the condition~$F\second\in\SS$ implies that
  $\smooth{F}\second\in\SS\subset\LL$ (and that in particular
  $\smooth{F}\second$ is tail-equivalent to $F\second$), we
  conclude that the
  supremum~$\smooth{M}$ of the process~$(\smooth{S}_n,\,n\ge0)$ satisfies
  \begin{equation}
    \label{eq:37}
    \limsup_{y\to\infty}\frac{P(\smooth{M}>y)}{\oo{F\second}(y)}
    \le \frac{C}{a-\epsilon}.
  \end{equation}
  It also follows from \eqref{eq:36} that the family of distributions
  $(\smooth{F}_x\,\,x\in\XX)$ satisfies the conditions~\eqref{eq:5}
  and \eqref{eq:6} of Proposition~\ref{iceland} with $F$ replaced by
  $\smooth{F}$, $\alpha$ replaced by $\delta$, and $\beta$ by
  $\beta+b$.  Hence, again since
  $\smooth{F}\second\in\SS$ and is tail-equivalent to $F\second$,
  there exists a constant~$r$ such that,
  for all $x\in\XX$, and for all $y$,
  \begin{equation}
    \label{eq:38}
    P(\smooth{M}>y\,\mid\,X=x) \le \min\left(1,r\oo{F\second}(y)\right).
  \end{equation}

  Define also the process~$(S^b_n,\,n\ge0)$ by 
  $S^b_n=\sum_{i=1}^n{}b_{X_i}$ for $n\ge0$.  Let $\eta_0=S^b_{T_0}$
  and $\eta_k=S^b_{T_k}-S^b_{T_{k-1}}$, $k\ge1$, be the increments of
  this process between the successive regeneration epochs of the
  modulating process~$X$.  It follows from \eqref{eq:33} that, for
  each $k\ge0$, $\eta_k\le{}b\tau_k$.  For a constant~$K$ to be
  specified below, define $\zeta_k = \max(\eta_k,b\tau_k-K)$ for each
  $k\ge0$.  The random variables $\zeta_k$ are independent for
  $k\ge0$, and are identically distributed for $k\ge1$; let $K$ be
  such that
  \begin{equation}
    \label{eq:39}
    E \zeta_k < 0, \qquad k \ge 1.
  \end{equation}
  Note also that $\zeta_k\le{}b\tau_k$ for each $k\ge0$ and so, from
  the condition~(D4),
  \begin{equation}
    \label{eq:40}
    P(\zeta_0 > y)  = o(\oo{F\second}(y)), \quad P(\zeta_1>y) = o(\oo{F}(y)),
    \qquad \text{as $y\to\infty$.}
  \end{equation}
  Now let $M^b=\sup_{n\ge0}S^b_n$.  Then
  \begin{align}
    M^b & \le \sup(b\tau_0,\eta_0+b\tau_1,\eta_0+\eta_1+b\tau_2,\dots)
    \nonumber\\
    & \le K +
    \sup(\zeta_0,\zeta_0+\zeta_1,\zeta_0+\zeta_1+\zeta_2,\dots)
    \nonumber\\
    & \le K + \zeta_0 + \sup(0,\zeta_1,\zeta_1+\zeta_2,\dots).
    \label{eq:41}
  \end{align}
  It follows from \eqref{eq:39}, \eqref{eq:40}, the independence of
  the random variables $\zeta_k$, $k\ge0$, and the Pakes-Veraverbeke Theorem
  that the last term on the right side of \eqref{eq:41} has a
  probability of exceeding $y$ which is $o(\oo{F\second}(y))$ as
  $y\to\infty$.  It now follows, from \eqref{eq:40}, \eqref{eq:41},
  the above independence and Lemma~\ref{addtails}, that
  \begin{equation}
    \label{eq:42}
    \lim_{y\to\infty}\frac{P(M^b>y)}{\oo{F\second}(y)} = 0.
  \end{equation}

  Finally, note that, for each $n$, we have $S_n=\smooth{S}_n+S^b_n$
  and hence $M\le\smooth{M}+M^b$. 
  Since $\smooth{M}$ and $M^b$ are conditionally independent
  given $X$, it follows from \eqref{eq:37}, \eqref{eq:38},
  \eqref{eq:42} and Lemma~\ref{addtails2} that
  \begin{equation}
    \label{eq:43}
    \limsup_{y\to\infty}\frac{P(M>y)}{\oo{F\second}(y)}
    \le
    \limsup_{y\to\infty}\frac{P(\smooth{M}+M^b>y)}{\oo{F\second}(y)}
    \le
    \frac{C}{a-\epsilon}.
  \end{equation}
  By letting $\epsilon\to0$ in \eqref{eq:43} and combining this result
  with the lower bound given by Lemma~\ref{lower}, we now obtain the
  required result.
\end{proof}

\subsection{Example}
\label{sec:example}

We give here an example to show the necessity of the condition (D4).

\begin{example}\label{example}
  Let $\zeta$, $\zeta_i$, $i\ge1$, be i.i.d.\ non-negative random
  variables with common distribution function $F$. Assume that
  $E\zeta=1$ and that $F\second\in\SS$.

  We take the modulating process $X=(X_n,\,n\ge 1)$ to be an
  independent Markov chain on $\Z_+=\{0,1,\dots\}$ with initial value
  $X_1=0$ and transition probabilities $p_{0,0}=0$, $p_{0,j}>0$ for
  all $j\ge1$ and, for $j\ge1$, $p_{j,j-1}=1$.  Define $T_0=0$ and for
  $k\ge1$, $T_k=\min\{n>T_{k-1}\colon{}X_n=1\}$.  We regard $T_k$,
  $k\ge0$, as the regeneration times of the process.  Since
  $p_{1,0}=1$, it follows that, for $k\ge1$, the $k$th cycle starts at
  time $T_{k-1}+1$ in state~$0$, and further that the cycle lengths
  $\tau_k=T_k-T_{k-1}$, are i.i.d.\ random variables with a
  distribution concentrated on $\{2,3,\dots\}$ and distribution
  function~$G$ given by $G(y)=\sum_{j\le{}y-1}p_{0,j}$.  Assume
  further that $E\tau_1=1+\sum_{j\ge1}jp_{0,j}<\infty$.  Then the Markov
  chain~$X$ is ergodic.

  Now define the modulated random walk $(S_n,\,n\ge0\}$ by $S_0=0$ and
  $S_n=\sum_{i=1}^n\xi_i$ where the random variables~$\xi_i$ are given
  by 
  \begin{displaymath}
    \xi_i = \zeta_i - d \1(X_i=0) 
  \end{displaymath}
  for some constant $d>E\tau_1$.  The conditions~(D1)--(D3) are thus
  satisfied with $F$ as defined here, $c(x)=1$ for all $x$, and
  $a=d/E\tau-1$.  
  
  Since the random variables~$\zeta_i$
  are nonnegative, we have $S_n\ge{}S_{n-1}$ for all $n$ such that
  $X_n\neq0$, i.e.\ for all $n\neq{}T_{k-1}+1$ for some $k$.  It
  follows that 
  \begin{equation}\label{eq:100}
    M \defn \sup_{n\ge0}S_n
    = \sup_{m\ge0} \left(\sum_{k=1}^m\psi_k \right)
  \end{equation}
  where, for $k\ge1$, 
  \begin{displaymath}
    \psi_k\defn\sum_{i=T_{k-1}+1}^{T_k}\zeta_i-d
  \end{displaymath}
  are i.i.d.\ random variables with common negative mean~$E\tau-d$.

  For the process~$(S_n,\,n\ge0)$ here, the constant~$\kappa$ defined
  by~\eqref{eq:2} is given $\kappa = E\zeta = 1$. For an arbitrary
  $b<1$, we provide an example when $P(b\tau>y)=o(\oo{F}(y))$, but for
  which $\oo{F\second}(y)=o(P(M>y))$, in each case as $y\to\infty$.
  Thus in this case the conclusion of Theorem~\ref{general} cannot
  hold.
  
  Choose $\gamma\in(0,1)$ and suppose that $\oo{F}(y)=e^{-y^{\gamma}}$
  (for which it is well-known that $F\second\in\SS$).  Suppose also
  that $P(\tau>y)\sim{}e^{-cy^{\gamma}}$ as $y\to\infty$ for some
  $c\in(b^{\gamma},1)$. Then it is readily checked that
  $P(b\tau>y)=o(\oo{F}(y))$.  For any fixed $\epsilon \in (0,1)$ such
  that $(1-\epsilon )^{\gamma}>c$, define the distribution~$H$ by
  \begin{displaymath}
    \oo{H}(y) = \exp \left(-\frac{cy^\gamma}{(1-\epsilon)^\gamma}\right).
  \end{displaymath}
  We now have
  \begin{align}
    P(\psi_1>y)
    & \ge
    P\left(
      \left| \frac{\sum_1^n \zeta_i}{n} -1 \right| \le \epsilon \quad
      \forall ~ n > y+d
    \right)
    P\left(\tau_1 > \frac{y+d}{1-\epsilon}\right)
    \nonumber\\
    & \sim P\left(\tau_1 > \frac{y+d}{1-\epsilon}\right)\label{eq:106}\\
    & \sim \oo{H}(y),\label{eq:107}
   \end{align}
   as $y\to\infty$, where \eqref{eq:106} follows by the Strong Law of
   Large numbers, and \eqref{eq:107} follows since $H\in\LL$.
   Since also $H\second\in\SS$, it follows from \eqref{eq:100},
   \eqref{eq:107}, and the Pakes-Veraverbeke Theorem (by for example noting
   that each random variable~$\psi_k$ stochastically dominates a random
   variable~$\psi'_k$ such that $P(\psi'_1>y)\sim \oo{H}(y)$) that
   \begin{displaymath}
     \liminf_{y\to\infty}\frac{P(M>y)}{\oo{H\second}(y)} \ge 1. 
   \end{displaymath}
   Finally, since also $\oo{F}(y)=o(\oo{H}(y))$, and so also
   $\oo{F\second}(x)=o(\oo{H\second}(y))$, as $y\to\infty$, it follows
   that $\oo{F\second}(x)=o(P(M>y))$ as required.

   Finally, we remark that while this example may be simplified
   somewhat by assuming the random variables~$\zeta_i$ to be a.s.\
   constant, we have some hope that, for a suitable choice
   of $F$, we may show the necessity of the \emph{strict}
   inequality~$b>\kappa$ in the condition~(D4).
\end{example}


\section{Modulated random walk in continuous time}
\label{sec:continuous-time}

In this section we consider a continuous-time process~$(S_t,\,t\ge0)$,
whose increments are independent and modulated by a background process
$X=(X_t,\,t\ge0)$ with a regenerative structure.  Analogously to the
discrete-time theory, the process is assumed to have jumps which are
heavy tailed and that $S_t\to-\infty$ as $t\to\infty$.  We are again
interested in the asymptotic form of the tail distribution of the
maximum of the process.  Many of the probabilistic ideas are similar
to the ones before.  However, we need to define the processes
carefully and we do so in Section~\ref{sec:deftriple}.  We then
present the main results in Section~\ref{sec:intr-main-results-1}, a
general result for processes with independent (but non-stationary)
increments in Section~\ref{sec:cts-iceland}, followed by the proofs in
Section~\ref{sec:proofs-theor-refc}.  We refer to Kallenberg (2002,
Ch.\ 15) for the theory and construction of processes with independent
increments.

\subsection{Definitions}
\label{sec:deftriple}

\paragraph{A process with independent increments.}
We define what we mean by a process
$(S_t,\,t\ge0)$ with $S_0=0$, independent increments and distribution
specified by a triple 
\[
(\nu(t,\cdot),\, \vr(t)^2,\,a(t),\quad t\ge0).
\]
First, $t \mapsto a(t)$ (respectively $t\mapsto \vr(t)^2$) is a
real-valued (respectively positive) function that is integrable over
finite intervals.  Second, for each $t$, the quantity $\nu(t,\cdot)$
is a Borel measure on $\R$ with $\nu(t,\{0\})=0$ and
$\int_\R(y^2\wedge|y|)\,\nu(t,dy)<\infty$; also, for each Borel set
$B$, the function $t \mapsto \nu(t, B)$ is integrable over finite
intervals.

Next let $\Phi$ be a Poisson random measure on $\R_+ \times \R$
with intensity measure $E\Phi(dt, dy) = dt~\nu(t, dy)$.
Note that the intensity measure is sigma-finite and so the Poisson 
random measure is well-defined.

Finally, for each $t$, let $A_t := \int_0^t a(s) ds$; let
$(W_t,\,t\ge0)$ be a zero-mean Gaussian process with independent
increments and $\var W_t = \int_0^t \vr(s)^2 ds$, and, for each
$t\ge0$, let $Y_t=\int_{[0,t]\times\R} y[\Phi(ds,dy)-ds\,\nu(s,dy)]$.
Note that the process~$(Y_t,\,t\ge0)$ is centred so that $EY_t=0$ for
all $t$.  Set
\begin{displaymath}
  S_t=A_t+W_t+Y_t, \qquad t\ge0.
\end{displaymath}
Thus, $S=(S_t,\,t\ge0)$ is a process with independent increments (see,
e.g., Kallenberg (2002)) and, in particular $ES_t=A_t$ for all $t$.  It is
not the most general version of a process with independent increments,
because we assumed that (i) its mean $ES_t =A_t$ exists (ii) the
functions $t \mapsto ES_t$, $t\mapsto E W_t^2$ are absolutely
continuous, and (iii) the intensity measure $E\Phi(dt,dy)$ has density
with respect to the first coordinate.  (Note that while the
assumptions~(ii) and (iii) are essentially technical, the
assumption~(i) is essential; in its absence we would need to pursue a
different treatment---in the spirit of Kl\"uppelberg, Kyprianou and
Maller (2004) and of Denisov, Foss and Korshunov (2004).)

\paragraph{A modulated continuous-time random walk.}
Next assume that we are given a regenerative process $X=(X_t,\, t \ge 0)$
such that $X_t$ takes values in some measurable space $(\XX, \XXX)$,
a measurable real-valued function $(a_x,\, x \in \XX)$, 
a measurable positive function $(\vr_x^2,\, x \in \XX)$,
and a collection of measures $(\nu_x(\cdot),\, x \in \XX)$, such
that $x \mapsto \nu_x(B)$ is measurable for each Borel set $B \subseteq \R$.

For each sample path $(X_t,\,t\ge0)$, define $S=(S_t,\, t \ge 0)$
as being a process generated by
the triple
\[
(\nu(t,\cdot),\, \vr(t)^2,\,a(t)) 
:= (\nu_{X_t}(\cdot),\, \vr_{X_t}^2,\, a_{X_t}).
\]
As above we assume that, for each $t$,
and each Borel set $B$,
\begin{equation}\label{eq:108}
  \int_0^t a_{X_s} ds, \quad
  \int_0^t \vr_{X_s}^2 ds, \quad
  \int_0^t \nu_{X_s}(B) ds
  \quad \text{are a.s.\ finite},
\end{equation}
and that, for each $x\in\XX$,
\[
\nu_x(\{0\})=0, \quad 
\int_\R (y^2 \wedge |y|) \nu_x(dy) < \infty.
\]
(For example we note that, if $\XXX$ is generated by some topology, a
sufficient condition for \eqref{eq:108} to hold is that $X$ have
\cadlag paths and that $x \mapsto a_x$, etc.\ be continuous functions.
)

We can construct this process by
considering a family $(\Phi_x,\,x \in \XX)$ of Poisson random measures
on $\R_+\times \R$ with intensity measure
\[
E \Phi_x(dt, dy) = dt ~\nu_x(dy),
\]
and an independent standard Brownian motion $(B_t,\,t\ge0)$.
We set
\begin{align*}
A_t &\defn \int_0^t a_{X_s}~ds \\
W_t &\defn \int_0^t \vr_{X_s}~dB_s \\
Y_t &\defn \int_0^t \int_\R y [\Phi_{X_s}(ds,dy) - \nu_{X_s}(dy)ds] \\
S_t &\defn A_t + W_t + Y_t.
\end{align*} 
We then have that, for each $t$, the characteristic function of $S_t$,
conditional on the background process $X$, is
\[
E [e^{i\theta S_t} \mid X]
=
\exp\left\{i\theta \int_0^t a_{X_s}~ds
-\frac{\theta^2}{2} \int_0^t \vr_{X_s}^2~ds
+\int_0^t~ds \int_\R~\nu_{X_s}(dy) [e^{i\theta y} -1 -i\theta y]\right\}.
\]
We shall refer to $S=(S_t,\,t\ge0)$ as a modulated continuous-time
random walk.  We assume that we choose a version of $S$ with \cadlag
paths.
(The reader will recognise that, in absence of modulation, the last
formula is the L\'evy-Khinchine formula for a L\'evy process--see
Bertoin (1998) or Sato (2000).)
We shall use the notation $\Delta S_t$ for the size of the
jump at any time $t$, i.e.\
\[
\Delta S_t \defn S_t-S_{t-}.
\]
We will also need to denote by $\Phi$ the point process 
on $\R_+ \times \R$ with atoms the pairs
$(t, \Delta S_t)$, for those $t$ for which $\Delta S_t \not = 0$, i.e.\
$\Phi(B) \defn \sum_{t: \Delta S_t \not = 0} \1((t,\Delta S_t)\in B)$,
$B \in \BB(\R_+\times \R)$.
Then, conditional on $X$, $\Phi$
is a Poisson point process with intensity measure
\begin{equation}
  \label{intensity}
  \BB(\R_+\times \R)
  \ni B \mapsto E [ \Phi(B) \mid X ]
  = \iint_B dt~ \nu_{X_t}(dy).
\end{equation}

\subsection{Main results}
\label{sec:intr-main-results-1}

We assume that the process~$(S_t,\,t\ge0)$ is modulated by a
regenerative background process~$X=(X_t,\,t\ge0)$ as described in the
previous section.

Denote the regeneration epochs of $X$ by $0\le{}T_0<T_1<\dots$. 
By definition, the cycles $((X_t,
T_{k-1} < t \le T_{k}),\, k\ge1)$ are i.i.d.\ and independent of the
initial cycle $(X_t,\, 0 < t \le T_0)$.  Define
$\tau_0=T_0$, $\tau_k=T_k-T_{k-1}$, $k\ge1$.
Then $(\tau_k,\,k\ge0)$ are independent and $(\tau_k,\,k\ge1)$ are
identically distributed.  We assume that $E\tau_1<\infty$.  For each
$t\ge0$, let $\pi_t$ be the distribution of $X_t$, 
and let $\pi$ denote the stationary probability measure
\begin{displaymath}
  \pi(B) \defn \frac{1}{E\tau_1} E \int_{T_1}^{T_2} \1(X_t \in B)~dt,\qquad
  B \in \BB({\mathcal X}).
\end{displaymath}

We require the extension of some definitions from distributions to
measures.  For any positive measure~$\nu$ on $\R$, again satisfying
the conditions
\begin{equation}
  \label{eq:109}
  \nu(\{0\})=0, \qquad \int_\R(y^2\wedge|y|)\,\nu(dy)<\infty, 
\end{equation}
we write $\oo{\nu}(y)\defn\nu((y,\infty))$ for all $y>0$.  We say that
$\nu$ is \emph{subexponential}, and write $\nu\in\SS$, if and only if
$\oo{\nu}(y)\sim{}c\oo{F}(y)$ as $y\to\infty$ for some
distribution~$F\in\SS$ and constant~$c>0$, i.e.\ if and only if $\nu$
has a subexponential tail; we similarly say that $\nu$ is
\emph{long-tailed}, and write $\nu\in\LL$, if and only if
$\oo{\nu}(y)\sim{}c\oo{F}(y)$ as $y\to\infty$ for some
distribution~$F\in\LL$ and constant~$c>0$.  Hence here also we have
$\SS\subset\LL$.  Finally, we define the \emph{integrated} (or
\emph{second}) \emph{tail} measure~$\nu\second$ on
$\R_+\setminus\{0\}$ by
$\oo{\nu\second}(y)=\int_y^\infty\oo{\nu}(z)\,dz<\infty$ for all
$y>0$.

The family $(\nu_x,\vr_x^2,a_x,\,x\in\XX)$ specifying the
distribution of $(S_t,\,t\ge0)$ will be assumed to satisfy the
following additional conditions with respect to some \emph{reference}
measure~$\nu$ on $\R$ satisfying~\eqref{eq:109}, some measurable
function $c:\XX\to[0,1]$ and constants~$\gamma$ and $\vr^2$:
\begin{flalign*}
  \text{\bf (C1)} & \quad
  \oo{\nu_x}(y)  \le \oo{\nu}(y)
  \quad \text{ for all $y>0$}, \qquad x\in\XX, &\\
  \text{\bf (C2)} & \quad
  \oo{\nu\second_x}(y)  \sim c(x) \oo{\nu\second}(y)
  \quad \text{ as $y \to \infty$}, \qquad x\in\XX, &\\
  \text{\bf (C3)} & \quad
  \int_{-\infty}^\infty (1\wedge y^2)\,\nu_x(dy) \le \gamma,
  \qquad x\in\XX, &\\
  \text{\bf (C4)} & \quad
  \vr_x^2 \le \vr^2, \qquad x\in\XX, &\\
  \text{\bf (C5)} & \quad
  \text{$a: = -\int_\XX a_x \pi(dx)$ is finite and strictly positive.}
\end{flalign*}

\begin{remark}\label{rem:conditionsC}
  The conditions~(C1), (C2) and (C5) are analogous to those of the
  discrete-time conditions~(D1), (D2) and (D3).  The remaining
  conditions~(C3) and (C4) are additional, and very natural,
  uniformity conditions necessitated by the continuous-time
  environment and have no (nontrivial) discrete-time analogues.  (With
  regard to the condition~(C3), note that the uniform boundedness in
  $x$ of $\oo{\nu_x}(1)$ is already guaranteed by the condition~(C1);
  the formulation of (C3) as above is for convenience.)
  Remarks~\ref{rem:gen}--\ref{rem:sign} in
  Section~\ref{sec:discrete-time} have obvious counterparts here.  We
  further remark that the condition imposed by (C3) on the left tails
  of the measures~$\nu_x$ may be weakened at the expense of some
  additional technical complexity.  Finally, note that only the
  restriction of the measure~$\nu$ to $\R_+\setminus\{0\}$ is relevant
  to the above conditions.
\end{remark}

As in Section~\ref{sec:discrete-time}, it follows from the above
conditions that the process~$(S_t,\,t\ge0)$ then satisfies
\begin{equation}
  \label{slln}
  \frac{S_t}{t} \to -a, \qquad\text{ as  $t\to\infty$, \quad a.s.}
\end{equation}
(see the discussion of this result in the Appendix).  Hence also
$S_t\to-\infty$, as $t\to\infty$, a.s., and so
\begin{math}
  M \defn \sup_{t\ge 0} S_t
\end{math}
is finite a.s. 

For each $x\in\XX$ and $\beta>0$, define
\begin{equation}
  \label{eq:89}
  a_x^\beta \defn a_x - \int_{-\infty}^{-\beta} (y+\beta)~\nu_x(dy).
\end{equation}
(Here $\int_{-\infty}^{-\beta}$ denotes $\int_{(-\infty,-\beta]}$; we
use similar conventions elsewhere.)  Define also
\begin{equation}
  \label{eq:44}
  \kappa\defn \lim_{\beta \to \infty} \sup_{x\in\XX} a_x^\beta.
\end{equation}
As in discrete time, in the case where the measures
$\nu_x$, $x\in\XX$, satisfy the uniform integrability condition
\begin{displaymath}
  \lim_{\beta \to \infty} \sup_{x\in\XX}
  \int_{-\infty}^{-\beta} |y|~\nu_x(dy) = 0,
\end{displaymath}
it follows from \eqref{eq:89} and \eqref{eq:44} that
$\kappa=\sup_{x\in\XX}a_x$.  Define $C\in[0,1]$ by
\begin{equation}
  \label{eq:45}
  C \defn \int_\XX c(x) \pi(dx).
\end{equation}
Theorem~\ref{cts-simple} below, for the case~$\kappa<0$, is the
analogue of Theorem~\ref{simple} for the discrete-time case.

\begin{theorem}
  \label{cts-simple}
  Suppose that (C1)--(C5) hold, that $\nu\second\in\SS$, and
  that $\kappa<0$.  Then
  \begin{displaymath}
    \lim_{y\to\infty}\frac{P(M>y)}{\oo{\nu\second}(y)}
    = \frac{C}{a}.
  \end{displaymath}
\end{theorem}

For the case where the sign of $\kappa$ may be arbitrary, we again
require an additional condition, similar to (D4), regarding the (tail)
distributions of the lengths of the regenerative cycles.  The
condition here is:

\text{\bf (C6)} \quad For some nonnegative $b>\kappa$,
\begin{equation}
  \label{eq:46}
  P(b\tau_0 > t)  = o(\oo{\nu\second}(t)),
  \quad P(b\tau_1>t) = o(\oo{\nu}(t)),
\qquad \text{as $t\to\infty$.}
\end{equation}
As in the discrete-time case, for $\kappa<0$ the condition~(C6) is
trivially satisfied by taking $b=0$, so that again
Theorem~\ref{cts-simple} may be viewed as a special case of the general
result given by Theorem~\ref{cts-general} below.

\begin{theorem}
  \label{cts-general}
  Suppose that (C1)--(C6) hold and that $\nu\second\in\SS$.  Then
  \[
  \lim_{y\to\infty} \frac{P(M>y)}{\oo{\nu\second}(y)} = \frac{C}{a}.
  \]
\end{theorem}


\subsection{A uniform upper bound for continuous-time processes with
  independent increments}
\label{sec:cts-iceland}

We prove in this section an auxiliary proposition, analogous to that
of Proposition~\ref{iceland} for the discrete-time case, which will be
required for the upper bound.  

\begin{proposition}
  \label{cts_iceland}
  Let $\nu$ be a Borel measure on $\R$ satisfying~\eqref{eq:109} and such
  that $\nu\second\in\SS$.  For strictly positive constants $\alpha$,
  $\beta$, $\gamma$, $\vr^2$, let the process~$(S_t,\,t\ge0)$ have
  distribution given by a triple $(\nu_t,\,\vr_t^2,\,a_t,\,t\ge0)$
  satisfying the conditions of Section~\ref{sec:deftriple} and such
  that, for all $t$,
  \begin{gather}
    \oo{\nu_t}(y) \le \oo{\nu}(y) \qquad\text{for all $y>0$}, \label{eq:47} \\
    \int_{-\infty}^{\infty} (1 \wedge y^2)\,\nu_t(dy) \le \gamma, \label{eq:49}\\
    \vr_t^2 \le \vr^2, \label{eq:50}\\
    a_t - \int_{-\infty}^{-\beta} (y+\beta)\,\nu_t(dy) \le -\alpha. \label{eq:51}
  \end{gather}
  Let $M\defn\sup_{t\ge0}S_t$.  Then there exists a constant $r$
  depending only on $\nu$, $\alpha$, $\beta$, $\gamma$ and $\vr^2$
  such that
  \begin{equation}
    \label{eq:52}
    P(M > y) \le r \oo{\nu\second}(y)
    \qquad\text{for all $y\ge0$}.
  \end{equation}
\end{proposition}

\begin{proof}
  Consider any process~$(S_t,\,t\ge0)$ with distribution given by
  $(\nu_t,\,\vr_t^2,\,a_t,\,t\ge0)$ as above.
  Choose $\epsilon\in(0,\alpha/2)$ and $y^*>0$ sufficiently large that
  \begin{equation}
    \label{eq:55}
    \oo{\nu}(y^*)\le\epsilon.
  \end{equation}
  Define, for each $t$, the measure
  $\nu^u_t$ by $\oo{\nu^u_t}(y)\defn\oo{\nu_t}(y^*\vee{}y)$---so that
  $\nu^u_t$ is the restriction of the measure~$\nu_t$ to
  $(y^*,\infty)$; define also, for each $t$, the (positive) measure
  $\nu^l_t$ by $\nu^l_t\defn\nu_t-\nu^u_t$---so that $\nu^l_t$ is the
  restriction of the measure~$\nu_t$ to $(-\infty,y^*]$.

  Decompose the process $(S_t,\,t\ge0)$ as $S_t=S^u_t+S^l_t$, where
  $S^u_0=S^l_0=0$, the process~$(S^u_t,\,t\ge0)$ has distribution
  given by $(\nu^u_t,\,0,\,-2\epsilon,\,t\ge0)$, and the
  process~$(S^l_t,\,t\ge0)$ is independent of $(S^u_t,\,t\ge0)$ and
  has distribution given by
  $(\nu^l_t,\,\vr_t^2,\,a_t+2\epsilon,\,t\ge0)$.  Define also
  $M^u\defn\sup_{t\ge0}S^u_t$ and $M^l\defn\sup_{t\ge0}S^l_t$.  Then
  $M^u$ and $M^l$ are independent and
  \begin{equation}
    \label{eq:56}
    M \le M^u + M^l.
  \end{equation}
  We now obtain upper bounds on the tail distributions of $M^u$ and
  $M^l$ which, in each case, depend only on $\nu$, $\alpha$, $\beta$,
  $\gamma$ and $\vr^2$.

  Define the measure $\bnd{\nu}$ concentrated on $(y^*,\infty)$ by
  $\oo{\bnd{\nu}}(y)\defn\oo{\nu}(y^*\vee{}y)$ for each $y>0$.  Since,
  for each $t$, $\nu^u_t$ is the restriction of the measure~$\nu_t$ to
  $(y^*,\infty)$ and similarly $\bnd{\nu}$ is the restriction of the
  measure~$\nu$ to $(y^*,\infty)$, it follows from \eqref{eq:47} and
  \eqref{eq:55} that, for all $y>0$, we have
  $\oo{\nu^u_t}(y)\le\oo{\bnd{\nu}}(y)\le\epsilon$.  Since also
  $ES^u_t=-2\epsilon{}t$ for all $t$, it follows (see
  Section~\ref{sec:deftriple}) that we may couple the process
  $(S^u_t,\,t\ge0)$ with a process~$(\bnd{S}_t,\,t\ge0)$, with
  $\bnd{S}_0=0$ and distribution given by the time-homogeneous triple
  $(\bnd{\nu},\,0,\,-\epsilon)$, in such a way that, almost surely,
  \begin{equation}
    \label{eq:57}
    S^u_t\le\bnd{S}_t \qquad\text{ for all $t$}.
  \end{equation}
  Define $\bnd{M}\defn\sup_{t\ge0}\bnd{S}_t$.  The process
  $(\bnd{S}_t,\,t\ge0)$ has i.i.d.\ positive jumps occurring as a
  Poisson process with
  rate~$\oo{\bnd{\nu}}(y^*)=\oo{\nu}(y^*)\le\epsilon$, and is linearly
  decreasing  between these jumps
  (i.e.\ it is a compound Poisson process with the subtraction of a
  linear function).  Let the random
  variables $0=t_0<t_1<t_2<\dots$ denote the successive jump times.
  Then the increments~$\bnd{\xi}_n=\bnd{S}_{t_n}-\bnd{S}_{t_{n-1}}$,
  $n\ge1$, of the process at the successive jump times are i.i.d.\
  random variables.  Since also, $E\bnd{S}_t=-\epsilon{}t$ for all
  $t$, we have $\oo{\nu}(y^*) E\bnd{\xi}_1\le-\epsilon$ and so
  $E\bnd{\xi}_1\le-1$.  Further the jumps of the process
  $(\bnd{S}_t,\,t\ge0)$ have a
  distribution~$G$ such that $\oo{G}(y)=\oo{\nu}(y^*\vee{}y)/\oo{\nu}(y^*)$.
  Since, as observed, the process is strictly decreasing between these
  jumps and
  since $\nu\second\in\SS\subset\LL$, it now follows from
  Lemma~\ref{lem:second} that the distribution~$H$ of $\bnd{\xi}_1$ is
  such that $\oo{H\second}(y)\sim\oo{G\second}(y)$ as $y\to\infty$.
  Hence, by the Pakes-Veraverbeke Theorem, there exists $\bnd{r}>0$ such that,
  for all $y\ge0$,
  \begin{equation}
    \label{eq:58}
    P(\bnd{M}>y) \le \bnd{r}\oo{\nu\second}(y).
  \end{equation}

  We now consider the tail distribution of $M^l$, and show that this
  is exponentially bounded.
  We show how to choose $s>0$, depending only on $\nu$, $\alpha$,
  $\beta$, $\gamma$ and $\vr^2$, such that the process
  $\left(e^{sS^l_t},\,t\ge0\right)$ is a supermartingale.  For this we
  require (from the distribution of $(S^l_t,\,t\ge0)$) that, for all $t$,
  \begin{equation}
    \label{eq:59}
    \frac{1}{s}\int_{-\infty}^\infty\left(e^{sy}-1-sy\right)\nu^l_t(dy)
    + a_t + 2\epsilon + \frac{\vr_t^2}{2}s \le 0.
  \end{equation}
  Define $K\defn\max(y^*,\beta,1)$.  We now use the upper bound, valid
  for any $s>0$,
  \begin{displaymath}
    \frac{1}{s}(e^{sy}-1-sy) \le
    \begin{cases}
      \frac{1}{s}(e^{-s\beta}-1-sy) \le -y-\beta + \frac{\beta^2}{2}s,
      & \qquad y \le -\beta,\\
      s e^{sK} y^2, 
      & \qquad -\beta < y \le y^*.
    \end{cases}
  \end{displaymath}
  Since also $\nu^l_t$ is the restriction of the measure~$\nu_t$ to
  $(-\infty,y^*]$, it follows that the left side of \eqref{eq:59} is
  bounded above by
  \begin{equation}
    \label{eq:60}
    \begin{split}
     \int_{-\infty}^{-\beta}\left(-y-\beta+\frac{\beta^2}{2}s\right)\nu_t(dy)
    + se^{sK}\int_{-\beta}^{y^*} y^2\,\nu_t(dy)
    + a_t + 2\epsilon + \frac{\vr_t^2}{2}s\\
    \le  -\alpha + 2\epsilon +s\left(\frac{\beta^2}{2}\gamma
    + e^{sK}K^2\gamma  + \frac{\vr^2}{2}\right),
   \end{split}
  \end{equation}
  by \eqref{eq:49}--\eqref{eq:51}, since, in particular,
  $y^2\le{}K^2(1\wedge{}y^2)$ on the interval~$(-\beta,y^*]$.  Finally,
  since $2\epsilon<\alpha$, it follows that $s$ may be chosen
  sufficiently small (and dependent only on $\nu$, $\alpha$, $\beta$,
  $\gamma$ and $\vr^2$) that the right side of \eqref{eq:60} is
  negative.

  As in the proof of Proposition~\ref{iceland}, it now follows, by the
  usual argument involving the martingale maximal inequality, that,
  for $s$ as above and $y\ge0$,
  \begin{equation}
    \label{eq:61}
    P(M^l>y) \le e^{-sy}.
  \end{equation}
  Now let $\zeta$ be random variable, independent of all else, which
  has tail distribution $e^{-sy}$.  From~\eqref{eq:56}, and since
  $M^u$ and $M^l$ are independent and (by construction) $M^u\le\bnd{M}$
  a.s., it follows that, for $y\ge0$,
  \begin{equation}\label{eq:62}
    P(M > y)
    \le P(M^u + M^l > y)
    \le P(M^u + \zeta > y)
    \le P(\bnd{M} + \zeta > y). 
  \end{equation}
  Again, as in the proof of Proposition~\ref{iceland}, it follows from
  the independence of $\bnd{M}$ and $\zeta$, \eqref{eq:58} and
  \eqref{eq:61}, and Lemma~\ref{addtails} that there exists $r$,
  depending only on $\nu$, $\alpha$, $\beta$, $\gamma$ and $\vr^2$,
  such that, for all $y>0$,
  \begin{displaymath}
    P(\bnd{M}+\zeta > y) \le r\oo{\nu\second}(y),
  \end{displaymath}
  and the required result now follows on using \eqref{eq:62}.
\end{proof}

\subsection{Proofs of Theorems~\ref{cts-simple} and \ref{cts-general}} 
\label{sec:proofs-theor-refc}

The following Lemma is analogous to Lemma~\ref{auxlem}.  Its proof is
entirely similar and so will be omitted.

\begin{lemma}
  \label{auxlem-c}
  Suppose that $\nu\second\in\LL$ and that $d_1$, $d_2$ are constants
  such that $d_2>0$.  Then
  \begin{equation} \label{eq:63}
    \int_0^\infty dt\; \oo{\nu}(y+d_1+d_2t)
    \sim \frac{1}{d_2}  \oo{\nu\second}(y)
    \qquad\text{as $y\to\infty$.}
  \end{equation}
  The conditions~(C1) and (C2) further imply that
  \begin{equation} \label{eq:64}
    \int_0^\infty dt \int_\XX\pi_t(dx)\;\oo{\nu_x}(y+d_1+d_2t)
    \sim \frac{C}{d_2} \oo{\nu\second}(y)
    \qquad\text{as $y\to\infty$.}
  \end{equation}
\end{lemma}

The following lemma gives an asymptotic lower bound for $P(M>y)$.

\begin{lemma}
  \label{lowerbound}
  Suppose that (C1)--(C5) hold and that $\nu\second\in\LL$.  Then
  \begin{displaymath}
    \liminf_{y\to\infty}
    \frac{P(M > y)}{\oo{\nu\second}(y)} \ge \frac{C}{a}.
  \end{displaymath}
\end{lemma}

\begin{proof}
  The proof of this is similar to that of Lemma~\ref{lower}.  Given
  $\epsilon>0$, by the weak law of large numbers we may choose a
  constant~$l_0$ sufficiently large that if, for each~$t$, we define
  $l_t=l_0+(a+\epsilon)t$, then
  \begin{equation}
    \label{eq:65}
    P(S_t>-l_t) > 1-\epsilon.
  \end{equation}
  Recall that, for each $t\ge0$,
  $S_{t-}\defn\lim_{u\uparrow{}t}S_{u}$ and
  $\Delta{}S_t\defn{}S_t-S_{t-}$; define also
  $M_{t-}\defn\sup_{u<t}S_u$.
  For each fixed $y\ge0$, note that the events
  \begin{displaymath}
    A_t \defn \{M_{t-}\le y,\,S_{t-}>-l_t,\,\Delta S_t>y+l_t\}
  \end{displaymath}
  (defined for all $t\ge0$) are disjoint---since each
  $A_t\subseteq\{M_{t-}\le{}y,\,M_t>y\}$.  Also, for each $t$, we
  have $A_t\subseteq\{M>y\}$.  Further, conditional on the background
  process~$X$, for each $t$, the events $\{M_{t-}\le
  y,\,S_{t-}>-l_t\}$ and $\{\Delta{}S_t>y+l_t\}$ are independent.  It
  follows that, for $y\ge0$,
  \begin{equation}
    P(M>y)  \ge P\Biggl( \bigcup_{t \ge 0} A_t \Biggr)
     = \int_0^\infty E[\1(\{M_{t-}\le y,\,S_{t-}>-l_t\})
    \oo{\nu_{X_t}}(y+l_t)]\, dt.
    \label{eq:66}
  \end{equation}
  (To obtain this result, we condition on the first (and only) time
  $t$ such that $\1_{A_t}=1$, and also use the fact that, conditional
  on $X$, the intensity measure of the point process~$\Phi$ introduced
  in Section~\ref{sec:deftriple} is as given by~\eqref{intensity}.)

  Now use the inequality, $E\1_A Z \ge EZ - c P(A^c)$, true for a
  random variable $Z$ such that $|Z| \le c$, a.s., to estimate the
  integrand in \eqref{eq:66} as
  \begin{multline*}
    E[\1(\{M_{t-}\le y,\,S_{t-}>-l_t\}) \oo{\nu_{X_t}}(y+l_t)]\\
    \ge E [\oo{\nu_{X_t}}(y+l_t)]
     -P(\{M_{t-} > y\} \cup \{S_{t-}\le-l_t\}) \oo{\nu}(y+l_t).
  \end{multline*}
  From this, \eqref{eq:66} and \eqref{eq:65}, we have that, as
  $y\to\infty$,
  \begin{align}
    P(M > y) & \ge \int_0^\infty E [\oo{\nu_{X_t}}(y+l_t)]\,dt
    -[P(M > y) + \epsilon] \int_0^\infty \oo{\nu}(y+l_t)\,dt \nonumber\\
    & = (1+o(1))\frac{C}{a+\epsilon}\oo{\nu\second}(y) - (1+o(1))[P(M
    > y) +
    \epsilon]\frac{\oo{\nu\second}(y)}{a+\epsilon} \label{eq:67}\\
    & = (1+o(1))\frac{C-\epsilon}{a+\epsilon}\oo{\nu\second}(y),
    \label{eq:68}
  \end{align}
  where \eqref{eq:67} follows from Lemma \ref{auxlem-c}, and
  \eqref{eq:68} follows since $P(M > y)\to0$ as $y\to\infty$.  The
  required result now follows by letting $\epsilon$ tend to zero.
\end{proof}

We now derive an asymptotic upper bound for $P(M>y)$ in the
case~$\kappa<0$.  The proof is similar to that of Lemma~\ref{upper}.

\begin{lemma}
  \label{cts-upper}
  Suppose that (C1)--(C5) hold, that $\nu\second\in\SS$, and that
  $\kappa<0$.  Then
  \begin{displaymath}
    \limsup_{y\to\infty}\frac{P(M>y)}{\oo{\nu\second}(y)}
    \le \frac{C}{a}.
  \end{displaymath}
\end{lemma}

\begin{proof}
  For given (small) $\epsilon>0$, and (large) $u_0>0$, define the
  linear function
  \begin{equation}
    \label{eq:69}
    u_t \defn u_0-(a-\epsilon)t, \qquad t \ge 0.
  \end{equation}
  Define the stopping time
  \begin{equation}
    \label{eq:70}
    \sigma \defn \inf\{t \ge 0:~ S_t > u_t\}.
  \end{equation}
  Since $S_t/t\to-a$ a.s., it follows that (for fixed $\epsilon$),
  \begin{equation}
    \label{eq:71}
    P(\sigma<\infty) \to 0 \qquad\text{as $u_0\to\infty$}.
  \end{equation}
  With regard to random variables such as $S_\sigma$ and
  $M_\sigma\defn\max_{0\le{}t\le{}\sigma}S_t$ which are only defined on
  $\{\sigma<\infty\}$, we again make the convention that, for example,
  $P(M_\sigma>y)\defn{}P(M_\sigma>y,\,\sigma<\infty)$.
  
  We first derive an upper bound for the tail of $M_\sigma$. 
  It follows from \eqref{eq:69} and \eqref{eq:70} that
  $M_{\sigma-}\le{}u_0$ a.s.\ on $\{\sigma < \infty\}$, and further
  that, for $y\ge{}u_0$,
  \begin{equation}
    \label{Ms}
    P(M_\sigma>y) = P(S_\sigma>y) \nonumber \\
    \le P(\Delta S_\sigma > y-u_\sigma).
  \end{equation}
  Let $\Phi$ be the point process 
  whose conditional intensity measure is given by \eqref{intensity}.
  Define also
  \begin{equation}
    \label{W}
    W = \{(t,z)\in \R_+ \times \R:~ z > y-u_t \}
  \end{equation}
 

  Note that if $\sigma<\infty$ and $\Delta{}S_\sigma>y-u_\sigma$, then
  $(\sigma,\Delta{}S_\sigma)\in W$ and hence the point process $\Phi$
  has at least one point in the region $W$.  Combining this last
  observation with the estimate \eqref{Ms}, we obtain
  \begin{align}
    P(M_\sigma > y) &\le P(\Phi(W)>0) \nonumber \\
    & = E[ P( \Phi(W) > 0 \mid X ) ] \nonumber \\
    & \le E[ E(\Phi(W) \mid X ) ] \label{eq:48} \\
    & = E \int_0^\infty \oo{\nu_{X_t}}(y-u_t) \,dt \nonumber \\
    & = \int_0^\infty dt \int_\XX \pi_t(dx) \oo{\nu_x} (y-u_t)\\
    & \sim\frac{C}{a-\epsilon}\oo{\nu\second}(y), \label{eq:72}
  \end{align}
  as $y\to\infty$, where \eqref{eq:48} follows since $\Phi(W)$ is
  nonnegative integer-valued, and \eqref{eq:72} follows from Lemma
  \ref{auxlem-c}.
  
  Since $\kappa<0$ it follows from \eqref{eq:44} that we can choose any
  $\alpha\in(0,-\kappa)$ and then $\beta>0$ sufficiently large that
  \begin{equation}\label{eq:73}
    a_x^\beta \le - \alpha
    \qquad\text{for all $x\in\XX$}.
  \end{equation}
  On the set $\{\sigma<\infty\}$ define the
  process~$(S^\sigma_t,\,t\ge0)$ by
  $S^\sigma_t\defn{}S_{\sigma+t}-S_\sigma$; let
  $M^\sigma\defn\sup_{t\ge0}S^\sigma_t$.  Conditional on the
  background process~$X$ and any finite value of $\sigma$, the
  process~$(S^\sigma_t,\,t\ge0)$ has independent increments and is
  generated by the triple
  $(\nu_{X_{\sigma+t}},\,\vr_{X_{\sigma+t}}^2,\,a_{X_{\sigma+t}},\,t\ge0)$.
  Further, it follows from the conditions~(C1), (C3), (C4), and
  \eqref{eq:73}, that, again conditional on $X$ and $\sigma$, this
  triple satisfies the conditions \eqref{eq:47}--\eqref{eq:51} of
  Proposition~\ref{cts_iceland} (with $\nu$, $\alpha$, $\beta$ as
  defined here and $\gamma$, $\sigma^2$ as defined by (C3) and (C4)).
  It therefore follows from Proposition~\ref{cts_iceland} that there
  exists a constant~$r$, depending on $\nu$, $\alpha$, $\beta$,
  $\gamma$ and $\vr^2$ only, such that, for all $x$, all finite $t$
  and all $y\ge0$,
  \begin{equation}
    \label{eq:74}
    P(M^\sigma>y \mid X=x,\,\sigma=t) \le r \oo{\nu\second}(y);
  \end{equation}
  further, conditional on $X=x$ and $\sigma=t$, the random variables
  $M_\sigma$ and $M^\sigma$ are independent.

  For $y>u_0$, we have $M=S_\sigma+M^\sigma=M_\sigma+M^\sigma$.  We
  now argue exactly as in the proof of Lemma~\ref{upper}, starting
  from the introduction of the random variable~$\widetilde{M}$ and
  with $F$ replaced by $\nu$ throughout, to obtain the required result.
\end{proof}


\begin{proof}[Proof of Theorem~\ref{cts-simple}]
  This is now immediate from Lemmas~\ref{lowerbound} and
  \ref{cts-upper}.
\end{proof}

\begin{proof}[Proof of Theorem~\ref{cts-general}]
  This is very similar to,
  but slightly simpler than, the proof of Theorem~\ref{general}.  Let
  nonnegative $b>\kappa$ be such that the condition~(C6) holds.
  Choose
  \begin{equation}
    \label{eq:1032}
    \delta\in(0,\min(a,b-\kappa))
  \end{equation}
  and choose $\epsilon\in(0,a-\delta)$.  Note that, from the
  condition~(C6),
  \begin{displaymath}
    \int_\XX(a_x + \delta)\pi(dx)= -a + \delta < -\epsilon.
  \end{displaymath}
  It now follows from the definition~\eqref{eq:44} of
  $\kappa$, and since $b>\kappa$, that we may choose $\beta>0$
  sufficiently large that
  \begin{gather}
    \int_\XX(a_x^\beta + \delta)\pi(dx) < -\epsilon, \label{eq:91}\\
    a_x^\beta + \delta \le b \qquad\text{for all $x\in\XX$}.
    \label{eq:88}
  \end{gather}
  Hence (as for example in the proof of Theorem~\ref{general}) we may
  define a measurable function~$b_x$ on $\XX$ such that,
  \begin{gather}
    a_x^\beta + \delta \le b_x \le b,
    \qquad x\in\XX, \label{eq:1033}\\
    \int_\XX b_x \pi(dx) = -\epsilon. \label{eq:1034}
  \end{gather}

  Define now the processes~$(S^b_t,\,t\ge0)$ and
  $(\smooth{S}_t,\,t\ge0)$ by, for each $t$,
  \begin{equation}
    \label{eq:94}
    S^b_t=\int_0^t{}b_{X_t},
    \qquad\qquad
    \smooth{S}_t= S_t - S^b_t.
  \end{equation}
  Note that, conditional on the background process~$X$, the
  process~$(\smooth{S}_t,\,t\ge0)$ has independent increments
  and a distribution which is given by the triple
  $(\nu_{X_t},\,\vr_{X_t}^2,\,\smooth{a}_{X_t},\,t\ge0)$, where, for
  each $x$, we have $\smooth{a}_x=a_x-b_x$.  It follows
  from~\eqref{eq:1034} that the process~$(\smooth{S}_t,\,t\ge0)$
  satisfies the conditions~(C1)--(C5) with $a$ is replaced by
  $a-\epsilon$.  Further, from the definitions~\eqref{eq:89},
  \eqref{eq:44} and the first inequality in \eqref{eq:1033}, the
  constant~$\kappa$ associated this process is replaced by some
  $\smooth{\kappa}$ satisfying $\smooth{\kappa}\le-\delta$.  Since
  also $\nu\second\in\SS$, it follows from Lemma~~\ref{cts-upper} that
  the supremum~$\smooth{M}$ of the process~$(\smooth{S}_t,\,t\ge0)$
  satisfies
  \begin{equation}
    \label{eq:1037}
    \limsup_{y\to\infty}\frac{P(\smooth{M}>y)}{\oo{\nu\second}(y)}
    \le \frac{C}{a-\epsilon}.
  \end{equation}

  It also follows from the conditions~(C1)--(C5) and the first
  inequality in~\eqref{eq:1033} that the family
  $(\nu_{X_t},\,\vr_{X_t}^2,\,\smooth{a}_{X_t},\,t\ge0)$ satisfies the
  conditions~\eqref{eq:47}--\eqref{eq:51} of
  Proposition~\ref{cts_iceland} with $\alpha$ replaced by $\delta$.
  Hence there exists a constant~$r$ such that, for all $x\in\XX$, and
  for all $y$,
  \begin{equation}
    \label{eq:1038}
    P(\smooth{M}>y\,\mid\,X=x) \le \min\left(1,r\oo{F\second}(y)\right).
  \end{equation}

  Now consider the process~$(S^b_t,\,t\ge0)$.  Recall that the
  condition~(C6) corresponds to the discrete-time condition~(D4) with
  $F$ replaced by $\nu$.  Recall also that $(T_k,\,k\ge0)$ is the
  sequence of regeneration epochs of the modulating process~$X$.  By
  considering the discrete-time process~$(S^b_{T_k},\,k\ge0)$, it
  follows \emph{exactly} as in the proof of Theorem~\ref{general}
  that, under the condition~(D7), the supremum~$M^b$ of the
  process~$(S^b_t,\,t\ge0)$ satisfies
  \begin{equation}
    \label{eq:1042}
    \lim_{y\to\infty}\frac{P(M^b>y)}{\oo{\nu\second}(y)} = 0.
  \end{equation}

  Finally, since $M\le\smooth{M}+M^b$, and since $\smooth{M}$ and
  $M^b$ are conditionally independent given $X$, it follows from
  \eqref{eq:1037}, \eqref{eq:1038}, \eqref{eq:1042} and
  Lemma~\ref{addtails2} that
  \begin{equation}
    \label{eq:1043}
    \limsup_{y\to\infty}\frac{P(M>y)}{\oo{\nu\second}(y)}
    \le
    \limsup_{y\to\infty}\frac{P(\smooth{M}+M^b>y)}{\oo{\nu\second}(y)}
    \le
    \frac{C}{a-\epsilon}.
  \end{equation}
  By letting $\epsilon\to0$ in \eqref{eq:1043} and combining this
  result with the lower bound given by Lemma~\ref{lower}, we now
  obtain the required result.
\end{proof}

\appendix

\section{Appendix}

In this appendix we give various general results concerning the
addition of subexponential random variables.  We also justify the
generalisations of the Strong Law of Large Numbers given by
\eqref{eq:97} and \eqref{slln}.

Lemma~\ref{addtails} below encapsulates the principle of one big jump
for subexponential random variables.  The result~\eqref{eq:80} is
standard---see, e.g., Baccelli, Schlegel and Schmidt (1999), while the
immediately following result follows by standard coupling arguments.

\begin{lemma}
\label{addtails}
Suppose that $F\in\SS$.  Let $Y_1, \dots, Y_n$ be independent random
variables such that, for each $i=1,\dots,n$, there exists a constant
$c_i>0$ with
\begin{equation}\label{eq:79}
  P(Y_i>y)\sim{}c_i\oo{F}(y)
  \qquad\text{as $y\to\infty$}
\end{equation}
(where in the case $c_i=0$ this is taken to mean
$P(Y_i>y)=o(\oo{F}(y))$ as $y\to\infty$).  Then
\begin{equation}\label{eq:80}
  P(Y_1+\cdots+Y_n > y) \sim (c_1+\dots+c_n) \oo{F}(y)
  \qquad\text{as $y\to\infty$}.
\end{equation}
Further, if, in \eqref{eq:79}, ``$\sim$'' is replaced by
``$\lesssim$'' for each $i$, then \eqref{eq:80} continues to hold with
``$\sim$'' similarly replaced by ``$\lesssim$''.
\end{lemma}

The following lemma gives a version of Lemma~\ref{addtails} (for the
case $n=2$ and with ``$\lesssim$'') where the random variables~$Y_1$
and $Y_2$ are conditionally independent.  It requires an extra,
asymmetric, condition (which is automatically satisfied in the case of
unconditional independence).

\begin{lemma}
  \label{addtails2}
  Suppose that $F\in\SS$.  Let $Y_1$ and $Y_2$ be random variables
  which are conditionally independent with respect to some
  $\sigma$-algebra~$\mathcal{F}$ and are such that, for some constants
  $c_1\ge0$, $c_2\ge0$, and some $r>0$,
  \begin{align}
    P(Y_i > y) & \lesssim c_i \oo{F}(y)
    \qquad\text{as $y\to\infty$},
    \qquad i=1,2, \label{eq:81}\\
     P(Y_1 > y\,\mid\,\mathcal{F}) & \le r\oo{F}(y)
    \qquad\text{for all $y$ \quad a.s.}  \label{eq:82}
  \end{align}
  (with the case $c_i=0$ interpreted as in Lemma~\ref{addtails}).  Then
  \begin{displaymath}
    P(Y_1+Y_2 > y) \lesssim (c_1+c_2) \oo{F}(y)
    \qquad\text{as $y\to\infty$}.
  \end{displaymath}
\end{lemma}

\begin{proof}
  Let $Y'$ be a random variable which is independent of $Y_2$ and such
  that
  \begin{equation}\label{eq:83}
    P(Y'>y) = 1 \wedge r\oo{F}(y) \qquad\text{for all $y$.}
  \end{equation}
  Since $F\in\SS$ implies $F\in\LL$, we can choose a positive
  increasing function $h_y$ of $y$ such that $h_y\to\infty$ as
  $y\to\infty$, but the convergence is sufficiently slow that
  $\oo{F}(y-h_y)\sim\oo{F}(y)$ as $y\to\infty$ (see, for example, Foss
  and Zachary (2002)).  Then
  \begin{align*}
    P(Y'+Y_2 > y)
    & =  P(Y_2 \le h_y,\,Y'+Y_2 > y) + P(Y_2 > h_y,\,Y'+Y_2 > y)\\
    & \le P(Y' > y-h_y) + P(Y_2 > h_y,\,Y'+Y_2 > y)\\
    & \sim r\oo{F}(y) + P(Y_2 > h_y,\,Y'+Y_2 > y) \qquad\text{as
      $y\to\infty$},
  \end{align*}
  where the last line above follows from \eqref{eq:82} and the
  definition of $h_y$.  Hence, since also, from \eqref{eq:81},
  \eqref{eq:82} and Lemma~\ref{addtails},
  $P(Y'+Y_2>y)\lesssim(r+c_2)\oo{F}(y)$ as $y\to\infty$, it follows
  that
  \begin{equation}
    \label{eq:84}
    P(Y_2>h_y,\,Y'+Y_2 > y) \lesssim c_2\oo{F}(y)
    \qquad\text{as $y\to\infty$}.
  \end{equation}
  We now have
  \begin{align}
    P(Y_1+Y_2 > y) & = P(Y_2 \le h_y,\,Y_1 + Y_2 > y) + P(Y_2 >
    h_y,\,Y_1 + Y_2 > y)
    \nonumber\\
    & \le P(Y_1 > y-h_y) + P(Y_2 > h_y,\,Y_1 + Y_2 > y)\nonumber\\
    & \le P(Y_1 > y-h_y) + P(Y_2 > h_y,\,Y' + Y_2 > y) \label{eq:85}\\
    & \lesssim P(Y_1 > y-h_y) + c_2\oo{F}(y)
    \qquad\text{as $y\to\infty$} \label{eq:86}\\
    & \lesssim (c_1+c_2) \oo{F}(y) \qquad\text{as $y\to\infty$},
    \label{eq:87}
  \end{align}
  as required, where \eqref{eq:85} follows by conditioning on
  $\mathcal{F}$, \eqref{eq:86} follows from \eqref{eq:84}, and
  \eqref{eq:87} follows from \eqref{eq:81} and the definition of
  $h_y$.
\end{proof}

Lemma~\ref{lem:second} below is a variant of a well-known
result.
\begin{lemma}
  \label{lem:second}
  Let $Y_1$ and $Y_2$ be independent random variables with
  distribution functions~$F_1$ and $F_2$ respectively.  Suppose that
  $F_1\second\in\LL$ and that $Y_2\ge0$ a.s ($F_2(y)=0$ for $y<0$).
  Then the distribution function~$F$ of $Y=Y_1-Y_2$ satisfies
  \begin{equation}
    \label{eq:96}
    \oo{F\second}(y) \sim \oo{F_1\second}(y)
    \qquad\text{as $y\to\infty$}.
  \end{equation}
  In particular, $F\second\in\LL$.
\end{lemma}

\begin{proof}
  The result is well-known when $F_1\second$ and $F\second$ in the
  statement of the lemma are replaced by $F_1$ and $F$
  respectively---see, e.g., Baccelli, Schlegel, and Schmidt (1999),
  and the modifications required for the present variation are
  trivially checked.
\end{proof}

Finally, we prove the generalisations of the Strong Law of Large
Numbers given by \eqref{eq:97} and \eqref{slln}.

Consider first the discrete-time case of
Section~\ref{sec:discrete-time}.  In the case where the modulating
process~$X$ is stationary (and, by definition, regenerative) then
$(\xi_n,\,n\ge0)$ is a stationary regenerative sequence and
\eqref{eq:97} follows from Birkhoff's theorem (since the invariant
$\sigma$-algebra is here trivial). In the general case,
one can always define a coupling of the sequence $(\xi_n,\,n\ge0)$ and
of a stationary regenerative sequence $(\xi'_n,\,n\ge0)$, such that
\begin{displaymath}
  \xi_{T_1+m} = \xi'_{T^{'}+m}
  \quad \text{a.s.} \qquad \text{for all $m=1,2,\dots$}
\end{displaymath}
for some non-negative and a.s.\ finite integer-valued random
variable~$T'$---see, for example, Thorisson (2000, Chapter 10, Section
3.)  Therefore, on the event $\{ T_1<n\}$,
\begin{displaymath}
  S_n = S'_{T^{'}-T_1+n} - S_{T^{'}}' + S_{T_1}
\end{displaymath}
and, as $n\to\infty$,
\begin{displaymath}
  \frac{S_n}{n} =
  \frac{S'_{T^{'}-T_1+n}}{T'-T_1+n}
  \frac{T'-T_1+n}{n} + \frac{S_{T_1}-S'_{T^{'}}}{n}
  \to -a \quad \text{a.s.}
\end{displaymath}
since the events $\{ T_1<n\}$ increase in $n$ to an event of
probability $1$.

The continuous-time result~\eqref{slln} follows entirely similarly.

\vspace*{1cm}
\noindent
\hfill
\begin{minipage}[t]{10cm}
\small \sc
Authors' address:\\
Department of Actuarial Mathematics and Statistics\\
School of Mathematical Sciences\\
Heriot-Watt University\\
Edinburgh EH14 4AS, UK\\
E-mail: {\tt \{S.Foss,T.Konstantopoulos,S.Zachary\}@ma.hw.ac.uk}
\end{minipage}


\begin{thebibliography}{20}
\bibitem{DFK}
{\sc Denisov, D., Foss, S.\ and Korshunov, D.} (2004).
Tail asymptotics for the supremum of a random walk
when the mean is not finite.
\emph{Queueing Systems}, \textbf{46}, 15--33.
\end{thebibliography}
\end{document}